\documentclass[trsc,nonblindrev]{informs3_ssrn} 

\OneAndAHalfSpacedXI

\usepackage{multirow}
\usepackage{array}
\usepackage{fix-cm}
\usepackage{mathptmx}
\usepackage{comment}
\usepackage{bm}
\usepackage{longtable}
\usepackage{diagbox}
\usepackage{booktabs,floatrow}
\usepackage{float}
\floatstyle{plaintop}
\restylefloat{table}
\usepackage{caption}
\usepackage{subcaption} 
\usepackage{stmaryrd}
\usepackage{natbib}
 \bibpunct[, ]{(}{)}{,}{a}{}{,}%
 \def\bibfont{\small}%
\usepackage{algorithm}
\usepackage{algpseudocode}
\usepackage{listings}
\newcommand{\norm}[1]{\left \lVert #1 \right \rVert}
\usepackage{appendix}

\TheoremsNumberedThrough     

\EquationsNumberedThrough    

\begin{document}




\TITLE{Physics-aware Truck and Drone Delivery Planning Using Optimization \& Machine Learning}

\ARTICLEAUTHORS{%
\AUTHOR{Yineng Sun}
\AFF{Thayer School of Engineering, Dartmouth College, Hanover, NH, United States of America, \EMAIL{yineng.sun.th@dartmouth.edu}, \URL{}}
\AUTHOR{Armin Fügenschuh}
\AFF{Brandenburgische Technische Universität Cottbus-Senftenberg, Cottbus, Germany, \EMAIL{fuegensc@b-tu.de}, \URL{}}
\AUTHOR{Vikrant Vaze}
\AFF{Thayer School of Engineering, Dartmouth College, Hanover, NH, United States of America, \EMAIL{vikrant.s.vaze@dartmouth.edu}, \URL{https://engineering.dartmouth.edu/community/faculty/vikrant-vaze}}
} 

\ABSTRACT{%
Combining an energy-efficient drone with a high-capacity truck for last-mile package delivery can benefit operators and customers by reducing delivery times and environmental impact. However, directly integrating drone flight dynamics into the combinatorially hard truck route planning problem is challenging. Simplified models that ignore drone flight physics can lead to suboptimal delivery plans. We propose an integrated formulation for the joint problem of truck route and drone trajectory planning and a new end-to-end solution approach that combines optimization and machine learning to generate high-quality solutions in practical online runtimes. Our solution method trains neural network predictors based on offline solutions to the drone trajectory optimization problem instances to approximate drone flight times, and uses these approximations to optimize the overall truck-and-drone delivery plan by augmenting an existing order-first-split-second heuristic. Our method explicitly incorporates key kinematics and energy equations in drone trajectory optimization, and thereby outperforms state-of-the-art benchmarks that ignore drone flight physics. Extensive experimentation using synthetic datasets and real-world case studies shows that the integration of drone trajectories into package delivery planning substantially improves system performance in terms of tour duration and drone energy consumption. Our modeling and computational framework can help delivery planners achieve annual savings worth millions of dollars while also benefiting the environment. 
}%

\HISTORY{}

\maketitle

\section{Introduction}
Recent years have witnessed explosive growth of the e-commerce sector. E-commerce sales increased by 6. 1\% in Q1 of 2025 from Q1 of 2024, accounting for 16.2\% of all sales \citep{census2024}.
Optimizing e-commerce delivery logistics has become increasingly crucial for improving profits and customer experience. Given the large number of packages to be delivered to customers in a timely and cost-effective manner, small but consistent improvements in delivery strategies can have a significant bottom-line impact.

Last-mile deliveries are typically performed by trucks, each operating in a predefined region within urban, suburban, or rural areas. Recent technological developments in unmanned aerial vehicles—also called drones—have led to new last-mile delivery schemes that pair trucks with drones to assist in package delivery. To describe the optimization problems in this setup, \citet{murray2015flying} coined the name, the flying sidekick traveling salesman problem (FSTSP). Since then, many studies have focused on different FSTSP variants, often using slightly different names and problem specifications. In particular, \citet{agatz2018optimization} defined and studied the traveling salesman problem with drones (TSP-D), where, for a single truck and a single drone, the delivery assignments and vehicle routes are optimized jointly. Unlike the arc-based modeling approach commonly used in FSTSP, TSP-D divides the entire truck-and-drone tour into segments called operations and uses binary variables to decide whether to select each operation to create the tour.

Package delivery optimization studies typically estimate drone travel times using simplified flight trajectory models. Recent advances in drone trajectory modeling capture some of the relevant flight physics to calculate travel times and drone energy consumption more accurately. \citet{schmidt2023two} modeled drone flights to generate precise, realizable and efficient trajectories, flight times, and drone energy consumption. Incorporating physics-based modeling of drone trajectories within the TSP-D framework can lead to a joint model that is more accurate but also more complicated. Given the associated computational difficulty, it is not yet clear to what extent such an integration might provide practical benefits. This paper shows that the common simplifications of drone flight physics used in truck-and-drone delivery optimization studies produce operational plans that can be highly suboptimal, demonstrating the need for better integration of nonlinear drone trajectory models within the combinatorial truck-and-drone planning problem.

Truck-and-drone delivery planning is a hard combinatorial optimization problem on its own, even with simplified drone-physics assumptions. Therefore, direct incorporation of drone trajectory equations into the base truck-and-drone delivery planning model is computationally intractable for realistic-sized problem instances. This paper proposes a drone-physics-aware truck-and-drone delivery planning model and an end-to-end solution approach that combines optimization and supervised learning to generate high-quality solutions in limited computational time. Extensive computational experimentation confirms that our approach scales well to real-world problem instances and outperforms computational and practical benchmarks consistently and substantially. Moreover, our model also incorporates restricted airspace constraints and enables drone take-off and landing on a moving truck to facilitate additional energy savings while ensuring safety.

This paper makes three contributions. First, we propose an integrated model of truck-and-drone delivery planning that combines nonlinear drone-trajectory physics with the combinatorial complexity of truck-and-drone route planning. Second, we develop an end-to-end solution approach that solves the resulting large-scale optimization model, using a neural network predictor combined with an augmented order-first-split-second routing heuristic. Third, we implement our modeling and algorithmic framework through a series of real-world case studies to demonstrate its practicality, substantial benefits, and operational insights.

\section{Related Literature}
This paper integrates nonlinear nonconvex relationships that describe drone flight physics into a combinatorial optimization model for truck-and-drone delivery planning. Our problem is at the intersection of the literature on truck-and-drone delivery planning and drone-flight physics, and our solution approach uses machine learning to accelerate transport logistics planning. Section \ref{lit:logistics} summarizes the related literature on truck-and-drone delivery planning, while Section \ref{lit:physics} reviews the work on drone-flight physics modeling. Section \ref{lit:ml-opt} reviews the key concepts and advances in machine learning to optimize transport logistics. Due to the large volume of literature in each of these areas, our objective is not to provide a comprehensive review. Instead, we focus on identifying the salient themes and positioning our work in their context.

\subsection{Truck-and-Drone Delivery Planning}\label{lit:logistics}

We can organize the literature on truck-and-drone delivery planning in terms of the various system configurations being considered, the mathematical formulations used, and the solution approaches employed.

\textbf{System Configurations:} \citet{murray2015flying} introduced optimization for collaborative delivery with a single truck and a single drone, a configuration adopted by many researchers \citep{bouman2018dynamic, ha2018min, marinelli2018route, liu2019cooperative, poikonen2020mothership, de2020variable, luo2022last, jeong2023drone}.Other studies have considered different setups, including multiple drones with one truck \citep{kitjacharoenchai2019multiple, murray2020multiple, raj2020multiple, luo2021multi, salama2022collaborative, jeong2023drone} and multiple drones with multiple trucks \citep{kitjacharoenchai2020two, masmoudi2022vehicle, yu2024electric, gao2023scheduling, rave2023drone, kloster2023multiple}. \cite{leon2022multi} allowed drones to visit more than one delivery node per flight.

Despite these differences, the primary goal of most models is to generate routing plans that use truck(s) and drone(s) to visit a specified set of delivery locations to optimize certain objectives. Although most studies use delivery locations directly as nodes in their models, some others also include auxiliary nodes under different names, such as auxiliary or intermediate depot, microdepot, central distribution center or flying warehouse \citep{jeong2020humanitarian, salama2022collaborative, jeong2023drone, rave2023drone}.

Truck-and-drone models typically approximate travel distances and travel times for both trucks and drones using straight-line Euclidean distances. Some studies estimate travel times that vary depending on the conditions of the road network \citep{wang2022truck}, or incorporate traffic congestion and related uncertainty to formulate robust truck-and-drone delivery routes \citep{yang2023planning}. Additional practical considerations include the effect of wind on drone flight \citep{sorbelli2023wind} and delivery service times \citep{das2020synchronized, yang2023planning}. Most studies allow drones to take off or land on a truck only at a delivery node or an auxiliary node. To account for recent advances in drone technology, some studies allow drones to take off and land on a moving truck anywhere along the truck path \citep{li2022truck, thomas2023collaborative}. We focus on a setting with one truck and one drone with the drone allowed to visit one delivery node per flight. Like most previous studies, we model delivery locations directly as model nodes. Additionally, we also account for restricted airspaces and allow the drone to take off and land when the truck is in motion or at rest. 

Collaborative truck-and-drone systems have been proposed for a variety of applications, including medical emergencies \citep{lin2022discrete}, disaster relief \citep{jeong2020humanitarian}, entertainment broadcasting \citep{ayranci2016use}, and surveillance \citep{zeng2022nested}. However, they are most commonly associated with commercial package delivery, which is the main focus of our work. We refer the reader to the review by \citet{moshref2021applications} focusing on logistics applications and those by \citet{chung2020optimization} and \citet{madani2022hybrid} for a broader overview of truck-and-drone systems with more general classifications.

\textbf{Mathematical Formulations:} The most common objective functions include the minimization of the total cost or the total duration of the tour. Some studies aim to optimize customer waiting time \citep{moshref2020design, moshref2021comparative}, operational cost \citep{gao2023scheduling, rave2023drone}, total cost to operators and users \citep{choi2021comparison}, weighted travel distance by vehicle type \citep{amorosi2021coordinating}, or total travel cost and customer satisfaction \citep{das2020synchronized}. Many existing studies use an arc-based modeling approach, which builds truck-and-drone trajectories by choosing whether or not to include each potential arc to connect a pair of nodes. Several others have adopted the operation-based approach of \citet{agatz2018optimization}. We also use the operation-based modeling approach to minimize the total duration of the tour.

\textbf{Solution Approaches:}
Even with simplified drone physics, off-the-shelf solvers are typically unable to handle real-world truck-and-drone delivery planning problems. Therefore, many researchers have developed custom solution methods. These include exact approaches \citep{amorosi2021coordinating, roberti2021exact, freitas2023exact, yang2023planning} and several heuristics such as the genetic algorithm \citep{remer2019multi, das2020synchronized}, simulated annealing \citep{moshref2020design,salama2022collaborative}, metaheuristics \citep{amorosi2021coordinating}, neighborhood search \citep{freitas2023exact, rave2023drone}, local search \citep{wang2022truck}, column generation-based heuristics \citep{gao2023scheduling}, and multi-agent heuristics \citep{leon2022multi}. The exact methods do not scale well—for example, the one proposed by \citet{agatz2018optimization} is tractable for instances with up to 12 nodes, although more than 100 nodes are common in practical package delivery planning problems. To address this challenge, we propose an end-to-end solution approach that scales well with increasing problem size and outperforms all tested benchmarks, consistently and substantially.

\subsection{Drone Flight Physics Modeling and Trajectory Planning} \label{lit:physics}
Last-mile drone delivery has received considerable attention in recent years. \citet{garg2023drones} review these drone applications, while \citet{macrina2020drone} provide a broad overview of drone usage in other contexts. \citet{pasha2022drone} review different modeling strategies for drone applications. Trajectory planning models vary in terms of the degree to which they account for physics and technical details \citep{fugenschuh2021flight}. Some studies accurately model the physical properties of drones \citep{ladyzynska2012modeling}, while others simplify them to prioritize other optimization objectives and handle larger instances \citep{adbelhafiz2010vehicle, pohl2008multi}. Most of the existing work on drone modeling lies somewhere between these two ends \citep{borrelli2006milp, jun2003path, geiger2006optimal}. \citet{schmidt2023two} present a mission and flight planning model, with a bilevel time structure to generate a smooth and practical drone trajectory. By linearizing the nonlinear constraints in their mixed-integer nonlinear optimization model, they solve the problem tractably using commercial solvers.

\subsection{Machine Learning to Optimize Transport Logistics Planning} \label{lit:ml-opt}
For larger and more complex optimization problem instances, solution approaches in the planning of transport and logistics operations must balance tractability with solution quality. Several studies have integrated machine learning (ML) to accelerate optimization \citep{bengio2021machine}. ML can be integrated into an optimization algorithm in multiple ways. One approach is to fully replace the optimization problem using ML, where the ML model trained on example problems and solutions directly outputs the final solutions to new instances \citep{larsen2018predicting}. Another approach is to use ML to augment large-scale optimization to enable more tractable decision-making \citep{rashedi2025machine, hizir2025machine}. A third strategy embeds the ML model within the optimization algorithm, repeatedly querying it for local decisions that improve the overall solution process \citep{khalil2017learning}. This paper falls into the second category, where we use ML to augment or configure optimization to obtain high-quality solutions quickly.

\citet{bengio2021machine} categorize ML-based optimization techniques by their properties. For well-defined large-scale optimization problems that require runtime improvement, methods that mimic or approximate solution patterns are ideal \citep{gasse2019exact, reiter2024equivariant}. In contrast, if the problem is challenging because certain expert knowledge is difficult to formalize or obtain, policy-based approaches that provide deeper insight into optimal solution structures are more suitable \citep{marcos2014supervised}. This paper falls into the first group—we develop neural network predictors to approximate drone travel times and use these estimates to inform optimization via instance-specific inputs.

\section{Model}
\citet{agatz2018optimization} defined an ``operation'' as a combination of a start node, an end end, at most one node served by the drone, and sequence of a non-negative number of nodes visited by either the truck alone or the truck-and-drone jointly (with drone resting on truck). Their model selects the operations that collectively form the final tour. They proposed an exact solution approach that is tractable for instances with up to 12 nodes. \citet{kundu2022efficient} proposed a polynomial-time splitting algorithm combined with local search that produces high-quality solutions quickly and scales to larger instances of the \citet{agatz2018optimization} model.

Like many previous studies, \citet{agatz2018optimization} and \citet{kundu2022efficient} assume that drone travel times are readily available, although estimating drone travel times is complicated in practice. Moreover, their models do not provide actual drone trajectories nor guarantee that an operation can be performed within its estimated time. To address these limitations, we introduce a new two-level formulation that improves solution quality and ensures the feasibility of the resulting truck-and-drone operations. Our model also allows the drone to take off or land on the truck anywhere along its path, reinforcing the modeling compatibility with drone technology. This capability further reduces travel time and drone energy consumption. We expand the definition of an operation from \citet{agatz2018optimization} to allow the drone to take off or land on a moving truck. Our model also prohibits the drone from entering restricted airspaces, according to local regulations.

\subsection{Problem Setup and Modeling Assumptions} \label{sec:setup}

Our setup considers a system of one truck and one drone to deliver packages to customers while minimizing the duration of the tour. Both vehicles begin at a depot node with the drone resting on the truck and all packages loaded on the truck. Our model makes the following assumptions.

\begin{itemize}
    \item The drone visits at most one delivery location in each operation. An operation begins when at least one vehicle departs from the operation's start node and ends when both vehicles reach the operation's end node. If the drone arrives first, it hovers to wait for the truck. If the truck arrives first, it waits for the drone.
    \item While resting on the truck, the drone is powered off. The drone battery is sufficient for its assigned tasks due to the sufficient battery capacity, the availability of charging on the truck, or battery swapping.
    \item To avoid conflict with ground traffic, the drone flies at a fixed altitude consistent with local regulations. It begins flight at the altitude of the truck bed, ascends and cruises at the target altitude, eventually descends to ground level and reaches zero velocity to drop off the package, and similarly returns to the truck.
    \item The truck driver or drone can easily locate and retrieve packages to serve any delivery node. Each delivery location corresponds to a single package. In case of multiple packages or multiple customers sharing a location, we redefine them as an aggregated single package for each location for modeling convenience.
    \item The truck travels at known speeds on each segment of the road, but different segments can have different speeds. So, the velocity vector of the truck is known for every segment of the road network.
    \item During the tour, the truck driver can park the truck at the delivery locations and the drone can fly to serve other delivery locations. Additionally, the drone is allowed to take off and land on the truck at any point during the truck's journey, including when the truck is in motion.
    \item The service time spent at each delivery location is negligible and hence omitted from the tour duration.
    \item In case of multiple solutions minimizing the duration of the operation, the optimal solution with the least drone energy consumption is used. Such ties are common in practice, making this tiebreaker important.
    \item The service area assigned to a truck-and-drone pair, the topology of the road network, and restricted airspaces are known well before the day of the tour. However, exact delivery locations are available only minutes or hours before the tour. So, ML models can be trained offline days ahead of operations, but route plans must be generated within minutes of obtaining the exact delivery locations for a given day's tour.
    \item Before a tour starts, the planner assigns sequences of nodes to the truck and the drone, but the exact trajectory of each flight is adjustable until that operation starts. This places a tighter runtime requirement on overall tour planning, but allows finalizing drone trajectories even after the truck-and-drone tour has started.
\end{itemize}

\subsection{Truck-and-Drone Tour Planning Model (TaDTPM)} \label{sec:TaDTPM}
The overall truck-and-drone tour planning problem decomposes into two sets of decisions. The upper level assigns delivery nodes to each vehicle and determines the visit sequences. The lower level decides the drone trajectories for all flights. The two levels are interdependent because the drone travel times decided by the lower-level drone trajectory model affect the total tour duration objective function of the upper-level model.

{
\footnotesize
\begin{longtable}{c c m{12.4cm}}
\caption{\footnotesize Upper-level model notation.} \label{notations_truckdrone_main}\\
\toprule
 \textbf{Symbol} & \textbf{Type} & \textbf{Description} \\
\toprule
 \endfirsthead

 \multicolumn{3}{c}{Continuation of Table \ref{notations_truckdrone_main}}\\
 \toprule
 \textbf{Symbol} & \textbf{Type} & \textbf{Description} \\
 \toprule
 \endhead
 \toprule
 \endfoot
 \endlastfoot

$\mathcal{U}$  & Set  & Nodes, including the delivery nodes, and the depot where the tour starts and ends\\
$\mathcal{O}$  & Set  & Feasible operations\\
$\mathcal{O}^-(u) \subseteq \mathcal{O}$  & Set  & Feasible operations that start at node $u \in \mathcal{U}$\\
$\mathcal{O}^+(u) \subseteq \mathcal{O}$  & Set  & Feasible operations that end in node $u \in \mathcal{U}$\\
$\mathcal{O}(u) \subseteq \mathcal{O}$  & Set  & Feasible operations that contain node $u \in \mathcal{U}$\\
$\overline{\mathcal{O}}(\mathcal{S}) \subseteq \mathcal{O}$  & Set & Feasible operations that start at a node in $\mathcal{U} \setminus \mathcal{S}$ and end in a node in $\mathcal{S}\subseteq \mathcal{U}$.\\
\hline
$u_{0}$ & Par. & Depot node.\\
$n$ & Par. & Number of delivery nodes.\\
\hline
$x_{o}$ & Var. & Binary. Indicates whether operation $o \in \mathcal{O}$ is performed.\\
$y_{u}$ & Var. & Binary. Indicates whether at least one operation starts at node $u \in \mathcal{U}$.\\
$t_{o}$ & Var. & Cont. Total duration of operation $o \in \mathcal{O}$.\\
$g_o^T$ & Var. & Cont. Drone energy consumed in operation $o \in \mathcal{O}$.\\
\bottomrule
\end{longtable}
}

Table \ref{notations_truckdrone_main} lists the notation in the upper-level model, Model \eqref{eq:ulmodel}, which is an extension of the model of \citet{agatz2018optimization}. Model \eqref{eq:ulmodel} is a lexicographic optimization model, which primarily minimizes tour duration, and from all duration-minimizing solutions, it secondarily chooses a solution that minimizes drone energy consumption. Compared to the linear objective in \citet{agatz2018optimization}, our primary objective, Equation \eqref{eq:sec3p5_objective13}, is a product of two variables. \citet{agatz2018optimization} assume operation times (denoted by $t_0$) to be constants, which makes the objective linear, while they are decision variables in our Model \eqref{eq:ulmodel} creating nonlinearity. In addition to this primary objective of total duration (TD) minimization, our model breaks the ties by introducing a secondary objective of minimizing drone energy consumption (DEC), Equation \eqref{eq:sec3p5_objective_secondary}, which is also nonlinear for similar reasons. Section \ref{lower-level-model} describes how the lower-level model computes $t_0$ and $g_0^T$.

Our model is based on the idea of an operation, defined as a segment of the overall truck-and-drone tour. The complete tour is a sequence of operations. Each operation starts at a start node with the drone on the truck and ends in an end node, again with the drone on the truck. An operation includes zero or one drone-only delivery node, plus any sequence of a non-negative integer number of truck nodes (visited by the truck alone or by the truck with the drone resting on top of it). Unlike \citet{agatz2018optimization}, our definition of an operation allows the drone to take off and land on the truck, even while it is in motion between nodes. Constraints \eqref{eq:sec3p5_constraint14}-\eqref{eq:sec3p5_constraint19} are similar to those in \citet{agatz2018optimization}. 
Constraints \eqref{eq:sec3p5_constraint14} ensure that all nodes are served. Constraints \eqref{eq:sec3p5_constraint15} define $y_{u}$ so that if at least one operation ends at node $u$, then at least one operation must start there. Constraints \eqref{eq:sec3p5_constraint16}
enforce the flow balance of the operations at each node. Constraints \eqref{eq:sec3p5_constraint17} eliminate sub-tours. Constraints \eqref{eq:sec3p5_constraint18} and \eqref{eq:sec3p5_constraint19}, respectively, require that at least one operation ends and starts at the depot node. Finally, the constraints \eqref{eq:sec3p5_constraint20} and \eqref{eq:sec3p5_constraint21} define the domains of the decision variables.

\begin{subequations}
\begin{align}
Primary: \qquad Minimize \qquad \sum_{o \in \mathcal{O}} t_o \cdot x_o & \label{eq:sec3p5_objective13}\\
Secondary: \qquad Minimize \qquad \sum_{o \in \mathcal{O}} g_o^T \cdot x_o & \label{eq:sec3p5_objective_secondary}\\
\sum_{o \in \mathcal{O}(u)} x_o \geq 1 & \qquad \forall u \in \mathcal{U} \label{eq:sec3p5_constraint14}\\
\sum_{o \in \mathcal{O}^+(u)} x_o \leq n \cdot y_{u} & \qquad \forall u \in \mathcal{U} \label{eq:sec3p5_constraint15}\\
\sum_{o \in \mathcal{O}^+(u)} x_o = \sum_{o \in \mathcal{O}^-(u)} x_o & \qquad \forall u \in \mathcal{U} \label{eq:sec3p5_constraint16}\\
\sum_{o \in \overline{\mathcal{O}}(\mathcal{S})} x_o \geq y_{u} & \qquad \forall \mathcal{S} \subset \mathcal{U}\setminus \{u_0\}, u \in \mathcal{S} \label{eq:sec3p5_constraint17}\\
\sum_{o \in \mathcal{O}^+(u_0)} x_o \geq 1 & \label{eq:sec3p5_constraint18}\\
y_{u_0} = 1 & \label{eq:sec3p5_constraint19}\\
x_{o} \in \{0,1\}, t_o, g_o^T \in \mathbb{R}_{\geq 0} & \qquad \forall o \in \mathcal{O} \label{eq:sec3p5_constraint20}\\
y_{u} \in \{0,1\} & \qquad \forall u \in \mathcal{U} \label{eq:sec3p5_constraint21}
\end{align} \label{eq:ulmodel}
\end{subequations}

\subsection{Drone Trajectory Model}\label{lower-level-model}
The lower-level model—an extension of the model of \citet{schmidt2023two}—optimizes the trajectory of the drone by coordinating with the location and velocity of the truck at take-off and landing. For each operation, this model takes as input an ordered set of truck delivery nodes (including operation start and end nodes) and a single drone delivery node, and produces the drone's trajectory, travel time, and energy consumption. It synchronizes the drone's trajectory and the truck's path at take-off and landing to allow take-off and landing anywhere along the truck's path. We discretize the planning horizon into time steps called \textit{major time steps}, and each major time step further consists of $n_f$ minor time steps, creating a two-level time discretization. The position, velocity, and acceleration of the drone are modeled in a fine-grained manner at each minor time step to ensure a high degree of modeling accuracy. On the other hand, variables describing the status of the drone (airborne or not) and its relative position with respect to the delivery location are modeled more coarsely using the major time steps to ensure tractability. This two-level discretization is a major contributor to the overall tractability and accuracy of the model. The reader is referred to \cite{schmidt2023two} for more details on this trade-off.

\begin{table}[hbt]
\footnotesize
\centering
\captionsetup{position=top} 
\caption{\footnotesize Lower-level model set notation.}\label{table_section3p1_set}
\begin{center}
\begin{tabular}{p{65mm}  p{9mm} p{85mm} }
\toprule
\textbf{Symbol} &  \textbf{Index}  &  \textbf{Definition}  \\
\midrule
$\mathcal{L}$  & $i$  & Drone altitude bands\\
$\mathcal{V}$ & $j$ & Drone throttle bands\\
$\mathcal{Q}$ & $q$ & Restricted airspaces\\
$\mathcal{T}=\{0,n_f,\cdots,n_f \cdot (T-1)\}$ & $t$ & Major time steps\\
$\mathcal{T}^{-}=\{0,n_f,\cdots,n_f \cdot (T-2)\}$ & $t$ & Major time steps excluding the last one\\
$\mathcal{T}_f=\{0,1,\cdots,n_f \cdot T-1\}$ & $t$ & Minor time steps\\
$\mathcal{T}_f^-=\{0,1,\cdots,n_f\cdot T-2\}$ & $t$ & Minor time steps excluding the last one\\
\bottomrule
\end{tabular}
\end{center}
\end{table}

{
\footnotesize
\begin{longtable}{c m{14cm}}
\caption{\footnotesize Lower-level model notation for parameters.} \label{notations_par}\\
\toprule
 \textbf{Symbol} & \textbf{Description} \\
\toprule
 \endfirsthead

 \multicolumn{2}{c}{Continuation of Table \ref{notations_par}}\\
 \toprule
 \textbf{Symbol} & \textbf{Description} \\
 \toprule
 \endhead
 \toprule
 \endfoot
 \endlastfoot

$n_f$ & Number of minor time steps in a major time step.\\
$T$  & Number of major time steps in the overall planning horizon.\\
$\tilde{T}$ & Last minor time step, $\tilde{T} = n_f \cdot T - 1$.\\
$\Delta_f$ & Length of a minor time step.\\
$P_i$  & Coordinate $i\in\{x,y,z\}$ of the drone delivery location.\\
$\lambda_2$ & 2D linearization constant = 0.3363788020.\\
$\lambda_3$ & 3D linearization constant = 0.2980450507.\\
$\delta$ & Maximum distance of drone to delivery node within which the node is considered visited.\\
$\overline{h}$ & Maximum flight altitude of drone.\\
$\underline{h}$ & Minimum flight altitude of drone.\\
$h^{min}$ & Minimum altitude of an airborne drone.\\
$R^0_i$ & Starting location of the drone along coordinate $i\in\{x,y,z\}$.\\
$R^{\tilde{T}}_i$  & Ending location of the drone along coordinate $i\in\{x,y,z\}$.\\
$\Delta$ & Length of a major time step.\\
$M$ & A sufficiently large number.\\
$\overline{v}^{z+}$, $\overline{v}^{z-}$  & Maximum climb and descent rates, respectively, of the drone.\\
$\overline{v}$, $\overline{a}$ & Maximum velocity and maximum acceleration, respectively, of the drone.\\
$F$ & Maximum drone energy level.\\
$m_v$  & Minimum drone charge level.\\
$\xi$ & Surplus (relative to horizontal flight) energy consumption rate per unit time per unit climb rate.\\
$\eta_{ij}$ & Energy consumption of the drone per unit time in the altitude band $i\in \mathcal{L}$, throttle band $j\in \mathcal{V}$.\\
$\theta_{j}$ & Discretized drone velocity value in throttle band $j \in \mathcal{V}$.\\
$H_{i}$ & Upper limit on the drone altitude in the altitude band $i \in \mathcal{L}$, where the default $H_0$ is 0.\\
$N^Q_q$ & Number of halfspaces whose intersection defines the $q^{\text{th}}$ RAS for $q \in \mathcal{Q}$.\\
$c_{qni}$  & Constant multiplier of $n^{\text{th}}$ half-space of $q^{\text{th}}$ RAS, along $i^{\text{th}}$ coordinate, for $q\in\mathcal{Q}$, $n\in\{1,\cdots,N^{Q}_{q}\}$, $i\in\{x,y,z\}$.\\
$c^{RHS}_{qn}$ & Right hand side constant for $n^{\text{th}}$ half-space of $q^{\text{th}}$ RAS, for $q\in\mathcal{Q}$, $n\in\{1,\cdots,N^{Q}_{q}\}$, $i \in \{x,y,z\}$. \\
$t_o^{trT}$ & Truck travel time for operation $o$.\\
$Z^{TR}$ & Truck driving altitude. \\
$t_o^{*}$ & Optimized total duration of the truck-and-drone operation $o \in \mathcal{O}$.\\
$\llfloor t \rrfloor$ & Major time step corresponding to minor time step $t\in \mathcal{T}_f$.\\
$PT_{it}$ & Coordinate $i\in\{x,y\}$ of truck location at minor time step $t\in \mathcal{T}_f$.\\
$VT_{it}$ & Truck velocity along coordinate $i \in \{x,y\}$ and minor time steps $t\in \mathcal{T}_f$.\\
\bottomrule
\end{longtable}

\begin{longtable}{c m{14cm}}
\caption{\footnotesize Lower-level model notation for decision variables. Note: Drone position vector is relative to delivery location.} \label{notations_var}\\
\toprule
 \textbf{Symbol}  & \textbf{Description} \\
\toprule
 \endfirsthead

 \multicolumn{2}{c}{Continuation of Table \ref{notations_var}}\\
 \toprule
 \textbf{Symbol} & \textbf{Description} \\
 \toprule
 \endhead
 \toprule
 \endfoot
 \endlastfoot

$b_t$ & Binary. Whether the drone is flying at major time step $t \in \mathcal{T}$.\\
$b^+_t$ & Binary. Whether drone has started flying at or before the major time step $t \in \mathcal{T}$.\\
$b^-_t$ & Binary. Whether the drone is available to fly and has not yet finished its flight at major time step $t \in \mathcal{T}$.\\
$w_{it}$ & Cont. Component $i \in \{x,y,z\}$ of drone's position relative to delivery location at major time step $t \in \mathcal{T}$.\\
$w^A_{it}$ & Cont. Absolute value of drone position vector's component $i \in \{x,y,z\}$ at major time step $t \in \mathcal{T}$.\\
$w^\infty_t$ & Cont. Maximum absolute value of the drone's 3D position vector components at major time step $t \in \mathcal{T}$.\\
$w^{L2}_t$ & Cont. Euclidean norm of the 3D position vector of the drone at major time step $t \in \mathcal{T}$.\\
$w^{M}_t$ & Cont. Maximum absolute value of horizontal components of drone's 2D position vector at major step $t \in \mathcal{T}$.\\
$w^{LA}_t$ & Binary. Whether absolute value of the $y$ component of the drone's position vector is larger than or equal to that of the $x$ component at major time step $t \in \mathcal{T}$.\\
$w^{LB}_t$ & Binary. Whether the absolute value of the $z$ component of the drone's position vector at major time step $t \in \mathcal{T}$ is larger than or equal to those of both the $x$ and $y$ components.\\
$w^D_{it}$ & Binary. Whether $i^\text{th}$ component of drone's position vector is negative at major time step $t\in\mathcal{T}$ for $i \in \{x,y,z\}$.\\
$d_t$ & Binary. Whether the drone visits the delivery location at major time step $t\in \mathcal{T}$.\\
$s_{ijt}$ & Binary. Whether the drone is in altitude band $i\in\mathcal{L}$ and throttle band $j\in\mathcal{V}$ at major time step $t \in \mathcal{T}$.\\
$f_{qnt}$ & Binary. Whether the drone is on the safe side of hyperplane $n\in \{1,\cdots,N^{Q}_{q}\}$ of RAS $q\in\mathcal{Q}$ at major time step $t\in\mathcal{T}$. To be outside of the RAS, the drone needs to be on the safe side of at least one hyperplane of that RAS.\\
$v^{z+}_t$ & Cont. Climb rate of the drone at minor time step $t\in \mathcal{T}_f$. It is set to zero when the drone is descending.\\
$v^{z-}_t$ & Cont. Descent rate of the drone at minor time step $t\in \mathcal{T}_f$. It is set to zero when the drone is climbing.\\
$v_{it}$ & Cont. Horizontal coordinate $i \in \{x,y\}$ of the drone velocity at minor time step $t \in \mathcal{T}_f$.\\
$v^A_{it}$  & Cont. Absolute value of component $i \in \{x,y\}$ of drone velocity at minor time step $t \in \mathcal{T}_f$.\\
$v^\infty_t$ & Cont. Maximum absolute value of the horizontal components of drone's velocity at minor time step $t\in \mathcal{T}_f$. \\
$v^{L2}_t$  & Cont. $l_2$ norm of the horizontal velocity of the drone at minor time step $t\in \mathcal{T}_f$.\\
$v^M_t$ & Binary. Whether $y$ component is the largest component of drone's horizontal velocity at minor time step $t\in\mathcal{T}_f$.\\
$v^D_{it}$ & Binary. Whether the component $i \in \{x,y\}$ of drone velocity is negative at minor time step $t\in\mathcal{T}_f$.\\
$a_{it}$ & Cont. Horizontal acceleration of the drone at minor time step $t \in \mathcal{T}_f$ along the component $i \in \{x,y\}$.\\
$a^A_{it}$ & Cont. Absolute value of components $i \in \{x,y\}$ of horizontal acceleration of drone at minor time step $t \in \mathcal{T}_f$.\\
$a^\infty_t$ & Cont. Maximum absolute value of horizontal components of drone's acceleration at minor time step $t \in \mathcal{T}_f$.\\
$a^{L2}_t$ & Cont. $l_2$ norm of the horizontal acceleration of the drone at minor time step $t \in \mathcal{T}_f$.\\
$a^M_t$ & Binary. Whether $y$ component is the largest component of the drone's horizontal acceleration at $t\in\mathcal{T}_f$.\\
$a^D_{it}$ & Binary. Whether the components $i \in \{x,y\}$ of drone acceleration is negative at minor time step $t\in\mathcal{T}_f$.\\
$r_{it}$ & Cont. The position of the drone along coordinate $i\in \{x,y,z\}$ at minor time step $t\in \mathcal{T}_f$.\\
$g_t$  & Cont. The amount of energy consumption of the drone at minor time step $t \in \mathcal{T}_f$.\\
\bottomrule
\end{longtable}
}

\vspace{-6mm}

Tables \ref{table_section3p1_set}, \ref{notations_par} and \ref{notations_var} define the sets, parameters, and decision variables, respectively, for this lower-level model. We now present the lower-level model for an operation $o \in \mathcal{O}$ in detail. Objective \eqref{eq:sec3p2_objective} is to minimize the total duration of the operation, which is defined by constraints \eqref{eq:sec3p2_oepr_time2_new}-\eqref{eq:sec3p2_oepr_time4}. Constraints \eqref{eq:sec3p2_oepr_time2_new} force $t_o$ to be at least equal to the end of the last major time step during which the drone is in the air, while the constraint \eqref{eq:sec3p2_oepr_time4} ensures that $t_o$ is at least equal to the duration of the truck for operation. Constraints \eqref{eq:sec3p2_touchz1}-\eqref{eq:sec3p2_touchv1} match the location and horizontal velocity vectors of the truck and drone whenever the drone rests on the truck. Constraints \eqref{eq:sec3p2_flightd1}-\eqref{eq:sec3p2_flightd3} (constraints \eqref{eq:sec3p2_flightdv1}-\eqref{eq:sec3p2_flightdv3}) ensure that the location (velocity) of the drone is updated depending on whether the drone is resting on the truck or flying. If the drone is airborne during at least one of the major time steps $\llfloor t\rrfloor$ and $\llfloor t+1\rrfloor$, then it must follow standard kinematic laws. Constraints \eqref{eq:sec3p1_airborne_c} ensure that $b_t=1$ if and only if both $b^+_t$ and $b^-_t$ are 1. Constraints \eqref{eq:sec3p1_airborne_takeoff} and \eqref{eq:sec3p1_airborne_landback} allow the drone to have a single transition from non-airborne to airborne status and a single transition from airborne to non-airborne status, respectively.

\vspace{-10mm}

\begin{subequations}
    \begin{align}
min \quad t_o & \label{eq:sec3p2_objective}\\
(t + n_f) \cdot \Delta_f \cdot b_t \leq t_o & \quad \forall t \in \mathcal{T} \label{eq:sec3p2_oepr_time2_new}\\
t_o^{trT} \leq t_o & \label{eq:sec3p2_oepr_time4}
\end{align}
\label{eq:2}
\end{subequations}

\vspace{-16mm}

\begin{subequations}
    \begin{align}
Z^{TR}- M \cdot b_{\llfloor t\rrfloor} \leq r_{zt} \leq Z^{TR}+ M \cdot b_{\llfloor t\rrfloor} & \qquad \forall t\in\mathcal{T}_f\label{eq:sec3p2_touchz1}\\
PT_{it} - M\cdot b_{\llfloor t\rrfloor} \leq r_{it} \leq  PT_{it} + M\cdot b_{\llfloor t\rrfloor} & \qquad \forall i \in \{x,y\}, t\in\mathcal{T}_f\label{eq:sec3p2_touchr1}\\
VT_{it} - M\cdot b_{\llfloor t\rrfloor} \leq v_{it} \leq  VT_{it} + M\cdot b_{\llfloor t\rrfloor} & \qquad \forall i \in \{x,y\}, t\in\mathcal{T}_f\label{eq:sec3p2_touchv1}
\end{align}
\label{eq:3}
\end{subequations}

\vspace{-12mm}

\begin{subequations}
\begin{align}
r_{it} + \Delta_f v_{it} + \frac{\Delta_f^2}{2} a_{it} -M(1-b_{\llfloor t+1\rrfloor}) \leq r_{i,t+1} \leq r_{it} + \Delta_f v_{it} + \frac{\Delta_f^2}{2} a_{it} +M(1-b_{\llfloor t+1\rrfloor}) & \:\:\:\:\:\: \forall i \in \{x,y\}, t \in \mathcal{T}^-_{f} \label{eq:sec3p2_flightd1}\\
r_{it} + \Delta_f v_{it} + \frac{\Delta_f^2}{2} a_{it} -M(1-b_{\llfloor t\rrfloor}) \leq r_{i,t+1} \leq r_{it} + \Delta v_{it} + \frac{\Delta_f^2}{2} a_{it} +M(1-b_{\llfloor t\rrfloor}) &\:\:\:\:\:\: \forall i \in \{x,y\}, t \in \mathcal{T}^-_{f} \label{eq:sec3p2_flightd3}\\
v_{it} + \Delta_f a_{it}-M(1-b_{\llfloor t+1\rrfloor}) \leq v_{i,t+1} \leq v_{it} + \Delta_f a_{it}+M(1-b_{\llfloor t+1\rrfloor}) &\:\:\:\:\:\: \forall i \in \{x,y\}, t \in \mathcal{T}^-_{f} \label{eq:sec3p2_flightdv1}\\
v_{it} + \Delta_f a_{it}-M(1-b_{\llfloor t\rrfloor}) \leq v_{i,t+1} \leq v_{it} + \Delta_f a_{it}+M(1-b_{\llfloor t\rrfloor}) &\:\:\:\:\:\: \forall i \in \{x,y\}, t \in \mathcal{T}^-_f \label{eq:sec3p2_flightdv3}
\end{align}
\end{subequations}

\vspace{-8mm}

\begin{subequations}
\begin{align}
b_t = b^+_t + b^-_t - 1 & \qquad \forall ~t\in\mathcal{T} \label{eq:sec3p1_airborne_c}\\
b^+_t \le b^+_{t+1} & \qquad  \forall t \in \mathcal{T}^- \label{eq:sec3p1_airborne_takeoff}\\
b^-_{t+1} \le b^-_t & \qquad \forall t \in \mathcal{T}^- \label{eq:sec3p1_airborne_landback}
\end{align}
\end{subequations}

Constraints \eqref{eq:sec3p1_dl1} define the distance vector from the delivery node to the current location of the drone at each major time step. Constraints \eqref{eq:sec3p1_dl2} require that, when the drone \textit{visits} a delivery node, it is within a predefined distance from the delivery location. Constraint \eqref{eq:sec3p1_dl3} ensures that the delivery location is visited exactly once during the drone flight. Constraints \eqref{eq:sec3p1_flightarea1} establish the upper and lower altitude limits for the drone, while constraints \eqref{eq:sec3p1_flightarea3} require that, except when delivering a package or at the start/end of its flight, the drone must remain above the minimum altitude required for airborne operations. Constraints \eqref{eq:sec3p1_se2}-\eqref{eq:sec3p1_se3} ensure that the operation starts and ends at the correct locations. Constraints \eqref{eq:sec3p1_flightdz1}-\eqref{eq:sec3p1_flightdz3} set the location and velocity of the drone along the $z$ axis and enforce upper limits on the climb rate and descent rate. Constraints \eqref{eq:sec3p1_flightdvmax} and \eqref{eq:sec3p1_flightdamax} enforce the upper limits of the drone's velocity and acceleration and ensure that the velocity and acceleration are set to zero when the drone is not airborne. Constraints \eqref{eq:sec3p1_energy1} define the upper and lower limits of drone energy consumption to ensure that the state of charge stays in the allowable range. Constraint \eqref{eq:sec3p1_energy3} initializes the energy consumption of the drone at 0\% of its capacity at the beginning of the operation. Constraints \eqref{eq:sec3p1_energy4} update the energy consumption of the drone at each minor time step. Constraints \eqref{eq:sec3p1_energy5} require that the drone is airborne if and only if it belongs to exactly one combination of the altitude band and throttle band at each major time step. Constraints \eqref{eq:sec3p1_energy6} and \eqref{eq:sec3p1_energy9}, respectively, ensure that the correct throttle band and the correct altitude band are selected at each minor time step. Constraints \eqref{eq:sec3p1_RAS1}-\eqref{eq:sec3p1_RAS2} prevent the drone from entering restricted airspace (RAS). Note that $q^\text{th}$ RAS is defined by $\{\bm{r}\in \mathbb{R}^3: \forall n\in \{1,\cdots,N_q^Q\}, \sum_{i \in \{x, y, z\}} c_{qni} \; \cdot \; r_i \ < c^{RHS}_{qn}\}$.

\vspace{-8mm}

\begin{subequations}
\begin{align}
w_{it}=r_{it}-P_i & \qquad \forall i \in \{x,y,z\}, t \in \mathcal{T} \label{eq:sec3p1_dl1}\\
w^{L2}_t \le \delta + M \cdot (1 - d_t) & \qquad \forall t \in \mathcal{T} \label{eq:sec3p1_dl2}\\
\sum_{t \in \mathcal{T}} d_t = 1 & \label{eq:sec3p1_dl3}
\end{align}
\end{subequations}

\vspace{-8mm}

\begin{subequations}
\begin{align}
\underline{h} \cdot b_{\llfloor t\rrfloor} \leq r_{zt} \le \overline{h} \cdot b_{\llfloor t\rrfloor} & \qquad \forall t \in \mathcal{T}_f \label{eq:sec3p1_flightarea1}\\
r_{zt} \geq h^{min} \cdot b_{\llfloor t\rrfloor} - M \cdot d_{\llfloor t\rrfloor} & \qquad \forall t \in \mathcal{T}_f \setminus \{0\}\label{eq:sec3p1_flightarea3}
\end{align}
\end{subequations}

\vspace{-8mm}
\begin{subequations}
\begin{align}
r_{i0}=R^0_i &\qquad \forall i \in \{x, y, z\} \label{eq:sec3p1_se2}\\
r_{iT}\geq R^{\tilde{T}}_i & \qquad \forall i \in \{x,y,z\} \label{eq:sec3p1_se3}\\
r_{z,t+1} = r_{z,t} + \Delta_f \cdot (v^{z+}_t - v^{z-}_t) & \qquad \forall t \in \mathcal{T}_f\label{eq:sec3p1_flightdz1}\\
v^{z+}_t \leq \overline{v}^{z+} &\qquad \forall t \in \mathcal{T}_f \label{eq:sec3p1_flightdz2}\\
v^{z-}_t \leq \overline{v}^{z-} & \qquad \forall t \in \mathcal{T}_f \label{eq:sec3p1_flightdz3}\\
v^{L2}_t \le \overline{v} \cdot b_{\llfloor t\rrfloor} & \qquad \forall t \in \mathcal{T}_f \label{eq:sec3p1_flightdvmax}\\
a^{L2}_t \le \overline{a} \cdot b_{\llfloor t\rrfloor} & \qquad \forall t \in \mathcal{T}_f \label{eq:sec3p1_flightdamax}
\end{align}
\end{subequations}

\vspace{-8mm}

\begin{subequations}
\begin{align}
0 \leq g_t \leq F-m_v & \qquad \forall t \in \mathcal{T}_f \label{eq:sec3p1_energy1}\\
g_0 = 0 & \label{eq:sec3p1_energy3}\\
g_{t+1} = g_t +\Delta_f( \xi * v^{z+}_t+\sum_{i\in\mathcal{L},j\in\mathcal{V}} \eta_{ij} s_{ij\llfloor t\rrfloor})  & \qquad \forall t\in\mathcal{T}^{-}_{f} \label{eq:sec3p1_energy4}\\
\sum_{i \in \mathcal{L},j \in \mathcal{V}} s_{ijt} = b_t & \qquad \forall t\in\mathcal{T} \label{eq:sec3p1_energy5}\\
v^{L2}_t -M\cdot(1-b_{\llfloor t\rrfloor}) \leq \sum_{i\in\mathcal{L},j\in\mathcal{V}} \theta_j s_{ij\llfloor t\rrfloor} \leq v^{L2}_t+M\cdot(1-b_{\llfloor t\rrfloor}) & \qquad \forall t\in\mathcal{T}_f \label{eq:sec3p1_energy6}\\
\sum_{i\in\mathcal{L}} H_{i-1} \sum_{j\in\mathcal{V}} s_{ij\llfloor t \rrfloor} - M\cdot(1-b_{\llfloor t\rrfloor})\leq r_{zt} \leq \sum_{i\in\mathcal{L}} H_i \sum_{j\in\mathcal{V}} s_{ij\llfloor t \rrfloor} +M\cdot(1-b_{\llfloor t\rrfloor}) & \qquad \forall t\in\mathcal{T}_f \label{eq:sec3p1_energy9}
\end{align}
\end{subequations}

\vspace{-8mm}

\begin{subequations}
\begin{align}
\sum_{i \in \{x, y, z\}} c_{qni} \; \cdot \; r_{it} \  \geq c^{RHS}_{qn} - M\left(1-f_{qn\llfloor t\rrfloor}\right) & \qquad \forall  q\in\mathcal{Q}, n\in \{1,\cdots,N^Q_q\}, t\in\mathcal{T}_f \label{eq:sec3p1_RAS1}\\
\sum_{n=1}^{N^Q_q} f_{qnt} \geq 1  & \qquad \forall  q\in\mathcal{Q}, t\in\mathcal{T} \label{eq:sec3p1_RAS2}
\end{align}
\end{subequations}

\textbf{Linearization of the Euclidean Norms}: Following \citet{schmidt2023two}, we approximate Euclidean ($l_2$) norm by a convex combination of the $l_1$ and $l_\infty$ norms \citep{chaudhuri1992modified} as $\norm{x}_2=\lambda_2\norm{x}_1+(1-\lambda_2)\norm{x}_\infty$ in two dimensions and as $\norm{x}_2=\lambda_3\norm{x}_1+(1-\lambda_3)\norm{x}_\infty$ in three dimensions, where $\norm{x}_1$, $\norm{x}_2$ and $\norm{x}_\infty$ represent the $l_1$, $l_2$, and $l_\infty$ norms, respectively, of a vector $x$. \cite{rhodes1995metrics} estimated the values as $\lambda_2 \approx 0.3363788020$ and $\lambda_3 \approx 0.2980450507$. The $l_1$ and $l_\infty$ norms themselves  are nonlinear, but they can be linearized with further reformulation by defining and using new binary variables. If $v_{it}\: \forall i \in \{x,y\}$ are the horizontal components of drone velocity at time $t$, the $l_2$ norm of the horizontal velocity vector is approximated by $\sqrt{v_{xt}^2+v_{yt}^2} \approx \lambda_2 \cdot \left(|v_{xt}|+|v_{yt}|\right) + (1-\lambda_2) \cdot max\{|v_{xt}|,|v_{yt}|\}$. We introduce the absolute value variables $v^A_{it}=|v_{it}|$, $\forall i \in \{x,y\}$, and denote the $l_2$ and $l_\infty$ norms as $v^{L2}_t=\sqrt{v_{xt}^2+v_{yt}^2}$, and $v^\infty_{t}=max\{|v_{xt}|,|v_{yt}|\}$. Binary variable $v^M_t=1$ if $v^\infty_{t}=v^A_{yt}$, and $v^M_t=0$ otherwise. Similarly, $v^D_{it}=0$ if $v^A_{it}=v_{it}$, and $v^D_{it}=1$ otherwise, $\forall i \in \{x,y\}$. Constraints \eqref{eq:sec3p1_normv1} approximate the $l_2$ norm as a linear combination of the $l_1$ and the $l_\infty$ norms, \eqref{eq:sec3p1_normv4}-\eqref{eq:sec3p1_normv7} set $v^\infty_{t}=max\{v^A_{xt},v^A_{yt}\}$, and \eqref{eq:sec3p1_normv2}-\eqref{eq:sec3p1_normv14} set $v^A_{it}=|v_{it}|, \forall i \in \{x,y\}$.

\vspace{-6mm}

\begin{subequations}
\begin{align}
v^{L2}_t = \lambda_2 \cdot \left(v^A_{xt} + v^A_{yt}\right) + (1-\lambda_2) \cdot v^\infty_t & \qquad \forall t \in \mathcal{T}_f \label{eq:sec3p1_normv1}\\
v^A_{it} \leq v^\infty_t &\qquad \forall i \in \{x,y\}, t \in \mathcal{T}_f \label{eq:sec3p1_normv4}\\
v^\infty_t - M \cdot v^M_t \leq v^A_{xt} \leq v^\infty_t + M \cdot v^M_t &\qquad \forall t \in \mathcal{T}_f \label{eq:sec3p1_normv5}\\
v^\infty_t - M \cdot (1-v^M_t) \leq v^A_{yt} \leq v^\infty_t + M \cdot (1-v^M_t) &\qquad \forall t \in \mathcal{T}_f \label{eq:sec3p1_normv7}\\
v_{it} \leq v^A_{it} & \qquad \forall i \in \{x,y\}, t \in \mathcal{T}_f \label{eq:sec3p1_normv2}\\
-v_{it} \leq v^A_{it} & \qquad \forall i \in \{x,y\}, t \in \mathcal{T}_f \label{eq:sec3p1_normv3}\\
v_{it} - M \cdot v^D_{it} \leq v^A_{it} \leq v_{it} + M \cdot v^D_{it} &\qquad \forall i \in \{x,y\}, t \in \mathcal{T}_f \label{eq:sec3p1_normv10}\\
-v_{it} - M \cdot (1- v^D_{it}) \leq v^A_{it} \leq -v_{it} + M \cdot (1-v^D_{it}) &\qquad \forall i \in \{x,y\}, t \in \mathcal{T}_f \label{eq:sec3p1_normv14}
\end{align}
\end{subequations}

\vspace{-2mm}

If $a_{it} \: \forall i \in \{x,y\}$ are the horizontal components of drone acceleration at time $t$, the $l_2$ norm of the horizontal acceleration vector is approximated by $\sqrt{a_{xt}^2+a_{yt}^2} \approx \lambda_2 \cdot \left(|a_{xt}|+|a_{yt}|\right) + (1-\lambda_2) \cdot max\{|a_{xt}|,|a_{yt}|\}$. We define continuous variables $a^A_{xt}, a^A_{yt}, a^{L2}_t, a^\infty_{t}$, and binary variables $a^M_t, a^D_{xt}, a^D_{yt}$ analogously to the corresponding variables for velocity. Constraints \eqref{eq:sec3p1_norma1}-\eqref{eq:sec3p1_norma14} linearize these acceleration relationships.

\vspace{-6mm}

\begin{subequations}
\begin{align}
a^{L2}_t = \lambda_2 \cdot \left(a^A_{xt} + a^A_{yt}\right) + (1-\lambda_2) \cdot a^\infty_t & \qquad \forall t \in \mathcal{T}_f \label{eq:sec3p1_norma1}\\
a^A_{it} \leq a^\infty_t &\qquad \forall i \in \{x,y\}, t \in \mathcal{T}_f \label{eq:sec3p1_norma4}\\
a^\infty_t - M \cdot a^M_t \leq a^A_{xt} \leq a^\infty_t + M \cdot a^M_t &\qquad \forall t \in \mathcal{T}_f \label{eq:sec3p1_norma5}\\
a^\infty_t - M \cdot (1-a^M_t) \leq a^A_{yt} \leq a^\infty_t + M \cdot (1-a^M_t) &\qquad \forall t \in \mathcal{T}_f \label{eq:sec3p1_norma7}\\
a_{it} \leq a^A_{it} & \qquad \forall i \in \{x,y\}, t \in \mathcal{T}_f \label{eq:sec3p1_norma2}\\
-a_{it} \leq a^A_{it} & \qquad \forall i \in \{x,y\}, t \in \mathcal{T}_f\label{eq:sec3p1_norma3}\\
a_{it} - M \cdot a^D_{it} \leq a^A_{it} \leq a_{it} + M \cdot a^D_{it} &\qquad \forall i \in \{x,y\}, t \in \mathcal{T}_f \label{eq:sec3p1_norma10}\\
-a_{it} - M \cdot (1-a^D_{it}) \leq a^A_{it} \leq -a_{it} + M \cdot (1-a^D_{it}) &\qquad \forall i \in \{x,y\}, t \in \mathcal{T}_f \label{eq:sec3p1_norma14}
\end{align}
\end{subequations}

\vspace{-2mm}

Finally, the $l_2$ norm of the drone's position vector at each major time step is approximated by constraints \eqref{eq:sec3p1_normw1}-\eqref{eq:sec3p1_normw23} using a strategy similar to that used for velocity and acceleration, except that the position vector is modeled in three dimensions as against the two-dimensional velocity and acceleration vectors modeled above. Another difference is that, as explained earlier, velocity and acceleration are modeled at each minor time step while the drone position relative to the delivery location is modeled at each major time step.

\vspace{-6mm}

\begin{subequations}
\begin{align}
w^{L2}_t = \lambda_3 \cdot (w^A_{xt}+w^A_{yt}+w^A_{zt}) + (1-\lambda_3) \cdot w^\infty_{t} & \qquad \forall t \in \mathcal{T} \label{eq:sec3p1_normw1}\\
w^{M}_t \leq w^\infty_{t} & \qquad \forall t \in \mathcal{T} \label{eq:sec3p1_normw6}\\
w^{A}_{zt} \leq w^\infty_{t} & \qquad \forall t \in \mathcal{T} \label{eq:sec3p1_normw7}\\
w^\infty_{t} - M \cdot w^{LB}_t \leq w^{M}_t \leq w^\infty_{t} + M \cdot w^{LB}_t & \qquad \forall t \in \mathcal{T} \label{eq:sec3p1_normw12}\\
w^\infty_{t} - M \cdot (1-w^{LB}_t) \leq w^{A}_{zt} \leq w^\infty_{t} + M \cdot (1- w^{LB}_t) & \qquad \forall t \in \mathcal{T} \label{eq:sec3p1_normw14}\\
w^{A}_{it} \leq w^{M}_t & \qquad \forall i \in \{x,y\}, t \in \mathcal{T} \label{eq:sec3p1_normw5}\\
w^{M}_t - M \cdot w^{LA}_t \leq w^{A}_{xt} \leq w^{M}_t + M \cdot w^{LA}_t & \qquad \forall t \in \mathcal{T} \label{eq:sec3p1_normw8}\\
w^{M}_t - M \cdot (1-w^{LA}_t) \leq w^{A}_{yt} \leq w^{M}_t + M \cdot (1- w^{LA}_t) & \qquad \forall t \in \mathcal{T} \label{eq:sec3p1_normw10}\\
w_{it} \leq w^{A}_{it} & \qquad \forall i \in \{x,y,z\}, t \in \mathcal{T} \label{eq:sec3p1_normw3}\\
-w_{it} \leq w^{A}_{it} & \qquad \forall i \in \{x,y,z\}, t \in \mathcal{T} \label{eq:sec3p1_normw4}\\
w_{it} - M \cdot w^D_{it} \leq w^{A}_{it} \leq w_{it} + M \cdot w^D_{it} &\qquad \forall i \in \{x,y,z\}, t \in \mathcal{T} \label{eq:sec3p1_normw17}\\
-w_{it} - M \cdot (1-w^D_{it}) \leq w^{A}_{it} \leq -w_{it} + M \cdot (1-w^D_{it}) &\qquad \forall i \in \{x,y,z\}, t \in \mathcal{T}\label{eq:sec3p1_normw23}
\end{align}
\label{eq:13}
\end{subequations}

\vspace{-6mm}

Finally, the constraints \eqref{eq:14domain01}-\eqref{eq:14domain12} define the domains of the decision variables.

\vspace{-6mm}
\begin{subequations}
\begin{align}
v^{z+}_t, v^{z-}_t, v^\infty_t, v^{L2}_t, a^\infty_t, a^{L2}_t,  g_t \geq 0 & \qquad \forall t \in \mathcal{T}_f\label{eq:14domain01}\\
w^\infty_t, w^{L2}_t, w^{M}_t \geq 0 & \qquad \forall t \in \mathcal{T}\label{eq:14domain02}\\
v_{it}, a_{it} \in \mathbb{R} & \qquad \forall i \in \{x,y\}, t \in \mathcal{T}_f\label{eq:14domain04}\\
v^A_{it}, a^A_{it} \geq 0 & \qquad \forall i \in \{x,y\}, t \in \mathcal{T}_f\label{eq:14domain04B}\\
w_{it} \in \mathbb{R} & \qquad \forall i \in \{x,y,z\}, t \in \mathcal{T}\label{eq:14domain05}\\
w^A_{it} \geq 0 & \qquad \forall i \in \{x,y,z\}, t \in \mathcal{T}\label{eq:14domain05B}\\
r_{it} \in \mathbb{R} & \qquad \forall i \in \{x,y,z\}, t \in \mathcal{T}_f\label{eq:14domain06}\\
v^M_t, a^M_t \in \{0,1\}  & \qquad \forall t \in \mathcal{T}_f\label{eq:14domain07}\\
b_t, b^+_t, b^-_t,  w^{LA}_t, w^{LB}_t, d_t \in \{0,1\}  & \qquad \forall t \in \mathcal{T}\label{eq:14domain08}\\
v^D_{it}, a^D_{it} \in \{0,1\}  & \qquad \forall i \in \{x,y\}, t \in \mathcal{T}_f\label{eq:14domain09}\\
w^D_{it} \in \{0,1\}  & \qquad \forall i \in \{x,y,z\}, t \in \mathcal{T}\label{eq:14domain10}\\
s_{ijt} \in \{0,1\}  & \qquad \forall i \in \mathcal{L}, j \in \mathcal{V}, t \in \mathcal{T}\label{eq:14domain11}\\
f_{qnt} \in \{0,1\}  & \qquad \forall q \in \mathcal{Q}, n \in \{1,\cdots,N^{Q}_{q}\}, t \in \mathcal{T}\label{eq:14domain12}
\end{align}
\label{eq:14_domain}
\end{subequations}

\vspace{-6mm}

In addition to the primary objective function of minimizing duration, we break ties by minimizing drone energy consumption (DEC) as a secondary objective function, modeled by \eqref{eq:sec3p3_objective}, and defined in constraint \eqref{eq:sec3p3_energy1}. Constraints \eqref{eq:sec3p3_oepr_time2_new} ensure that while minimizing DEC, the drone still respects the least duration established by the time-minimizing truck-and-drone model (Equations \eqref{eq:2}-\eqref{eq:14_domain}). Recall that $t_o^{*}$ is the optimized total duration of the operation obtained by solving the duration-minimizing truck-and-drone model. Overall, we first set $t_o^{*}$ at the lowest possible value of $t_o$ given by solving the model \eqref{eq:2}-\eqref{eq:14_domain}. Subsequently, we solve the tie-breaker model, which has the objective function given by the equation \eqref{eq:2} and the constraints given by equations \eqref{eq:3}-\eqref{eq:14_domain}, \eqref{eq:sec3p3_oepr_time2_new}, and \eqref{eq:sec3p3_energy1}.

\vspace{-10mm}

\begin{subequations}
    \begin{align}
min \quad g_o^T & \label{eq:sec3p3_objective}\\
(t + n_f) \cdot \Delta_f \cdot b_t \leq t_o^{*} & \quad \forall t \in \mathcal{T} \label{eq:sec3p3_oepr_time2_new}\\
g_o^T=g_{\tilde{T}} & \label{eq:sec3p3_energy1}
\end{align}
\end{subequations}

The total energy consumption of the drone for the entire truck-and-drone tour is the sum of the DECs of all selected operations. Formally, if $g_o^T$ is the optimal energy consumption of the drone for operation $o \in O$, then the total drone energy consumption for the entire tour is given by $\sum_{o \in O} g_o^T \cdot x_o$.

\section{Solution Approach}
Our two-level decomposition achieves two goals. First, for a given set of lower-level objective function coefficients (namely $t_o$ and $g_o^T, \forall o \in \mathcal{O}$), it ensures that the upper level is a biobjective mixed-integer linear optimization problem. Second, although the lower level is a biobjective mixed-integer nonconvex optimization problem, the decomposition substantially reduces its size: Without the decomposition, the overall optimization problem needs to be solved to jointly obtain the optimal trajectories for \textit{each} candidate operation, whereas the decomposition breaks it down into individual optimization problems, one for each candidate operation. For fixed values of $t_o$ and $g_o^T, \forall o\in\mathcal{O}$, the upper level is a biobjective extension of the mixed-integer linear optimization model of the TSP-D \citep{agatz2018optimization}.  \citet{kundu2022efficient} proposed a polynomial time splitting algorithm combined with local search for the single-objective TSP-D and produced high-quality solutions to realistic problem sizes. If the objective function coefficients are known beforehand, the TSP-D model can be solved by extending this approach. However, accurate estimation of drone travel times requires solving the lower-level model for all candidate operations, and the number of such operations grows exponentially with the number of delivery locations, making a full enumeration unrealistic.

Instead, we propose a new neural network-based approach to generate accurate estimates of drone travel time for each operation. \citet{kundu2022efficient} use Euclidean travel distances, which we propose to replace with a neural network (NN). Our NN model is trained on offline data and outputs a drone travel time estimate given the coordinates of three locations: operation start, drone delivery, and operation end. Next, we combine these NN-based estimates with an augmented version of the heuristic by \citet{kundu2022efficient} to generate the optimal set of operations. Finally, the detailed lower-level model can be used to provide a feasible 4D ($x$, $y$, $z$, and $t$) trajectory plan, duration, and DEC for each selected operation. We present our overall solution approach in Section \ref{sec:4.0}. Section \ref{sec:4.1} describes the NN design and training process. Section \ref{sec:4.2} presents the drone-only trajectory model used to generate the NN training data. Finally, Section \ref{sec:4.3} provides three methods for calculating drone flight times with varying degrees of simplicity.

\subsection{Overall Solution Approach} \label{sec:4.0}
The algorithm of \citet{agatz2018optimization} was improved by \citet{kundu2022efficient} (see \citet{kundu2022developing} for implementation details). The general concept behind these solution approaches has been summarized as ``order-first split-second'' or ``route and reassign''. Additionally, since the solution can be improved using an improvement heuristic as a third step, these approaches can be considered ``order-first, split-second, improve-third'' methods. We also use a ``order-first, split-second, improve-third'' approach, with the following three steps.

\textbf{Ordering}: In the first step, these approaches assume that the drone does not leave the truck throughout the tour, and instead the truck covers all delivery nodes as parts of a single giant tour. Because this is a TSP, known to be NP-hard, and because the node ordering in this step only indirectly determines the quality of the final truck-and-drone solution, heuristics are typically used to obtain a fast initial TSP solution. \citet{agatz2018optimization} used the Concorde TSP solver \citep{applegate2006traveling}, and \citet{kundu2022developing} used the Lin-Kernighan-Helsgaun (LKH) algorithm \citep{helsgaun2000effective}. We chose to use the two-opt local search heuristic \citep{croes1958method} for its shorter runtime. We also experimented with other TSP heuristics, but they neither improved the quality of the final truck-and-drone solution nor affected the main conclusions and takeaways.

\textbf{Splitting}: The splitting step evaluates—on the basis of their durations, obtained by the estimation methods in Section \ref{sec:4.3}—all possible operations in which the nodes visited by the truck within that operation are in the same order as that induced by the truck-only tour. For this splitting step, \citet{agatz2018optimization} used an exact dynamic programming-based partitioning method, while \citet{kundu2022efficient} used a shortest path procedure on a directed acyclic graph to reduce computational overhead without sacrificing optimality. Due to its computational efficiency, we chose the splitting method in \citet{kundu2022efficient}, whose key idea is to keep track of the minimum travel time to reach each node from the start node and update this minimum travel time during the evaluation process whenever new operations with a shorter travel time are encountered. This step provably produces the shortest possible truck-and-drone tour among all tours where the truck visits nodes in the same order as that induced by the truck-only tour from the first step.

Since there are $O(n^3)$ possible operations with the same induced ordering as a given $n$-node truck-only tour, the evaluation process requires at least $O(n^3)$ runtime. Depending on the method used to calculate the duration of each operation, the overall runtime can be $O(n^4)$ as in the exact partitioning method in \citet{agatz2018optimization} or $O(n^3)$ as in the speed-up version of the method in \citet{agatz2018optimization}, as well as in the method proposed in \citet{kundu2022efficient} and \citet{kundu2022developing}. We use the splitting step as implemented in \citet{kundu2022developing} where the operation time calculation takes $O(1)$ runtime, and thus the splitting step has $O(n^3)$ time complexity. We use this method because it gives higher quality solutions with runtimes of more than 95\% lower than the methods in \citet{agatz2018optimization} for large problem instances \citep{kundu2022efficient}.

\textbf{Improvement}:
Because the truck-and-drone tour obtained by the splitting step depends on the initial truck-only tour, it can be improved by considering alternative truck-only tours to initiate the splitting step. To find ``neighbors'' of the initial truck-only tour, both \citet{agatz2018optimization} and \citet{kundu2022efficient} consider an iterative improvement procedure involving three different methods: one-point move (relocate one node), two-point move (swap two nodes), and two-opt (replace two edges). We run the splitting step on all neighbors of the initial truck-only TSP tour to find their corresponding shortest truck-and-drone tours and pick the best one, which improves the truck-and-drone solution compared to the original single run of the splitting method. See Section 6.3 in \citet{agatz2018optimization} for more details on this improvement step.

Importantly, unlike \citet{agatz2018optimization} and \citet{kundu2022developing}, we do not estimate drone travel times as Euclidean distances divided by a constant speed. Sections \ref{sec:4.1}-\ref{sec:4.3} detail our estimation method.

\subsection{Neural Network Model Predictor for Drone-Only Time Estimation} \label{sec:4.1}
Computing a specific drone trajectory and travel time requires solving the lower-level model \eqref{eq:2}-\eqref{eq:14_domain}, a calculation repeated several times when solving large-scale problems, which is computationally demanding. We use a NN model to accelerate the runtimes. The input layer of the NN encodes the fundamental descriptors of a drone operation: the $x$ and $y$ coordinates of the drone's start, delivery, and end locations within an operation, comprising six features. The NN produces one output: drone travel time. To train the NN, we generate instances of operations assuming uniformly distributed locations of all three nodes, and run the drone-only model described in Section \ref{SA:DOSDTM} to obtain the ``true'' drone travel time labels for each instance. We then pair these labels with the instance's six-feature vector to build a training dataset for the NN.

We prefer a simple NN structure, with a single hidden layer to reduce the complexity and runtime of the training and tuning processes. A single hidden layer with a sufficient number of neurons can approximate any continuous function, and structures with a single hidden layer are common in regression \citep{hornik1989multilayer}. Moreover, function approximation tasks with high-dimensional data and a large number of neurons generalize well \citep{zhang2017understanding}. Therefore, we choose to use a neural network with a large single hidden layer. We tune the number of neurons in this hidden layer and other hyperparameters via a standard grid search. For each hyperparameter considered in the grid search, we first refer to the default values from the scikit-learn package called Multilayer Perceptron Regressor, or MLPRegressor, and fine-tune the hyperparameters around these default values. As too few neurons lead to underfitting, while too many can cause overfitting and unnecessary computational cost \citep{goodfellow2016deep}, we consider a range of choices (1000, 2500 or 4000 hidden neurons) to balance expressiveness and overfitting risks.

For the activation function, we provide options that are suitable for both near-linear cases (Identity) and nonlinear cases (ReLU) \citep{glorot2011deep} to offer flexibility. The alpha value controls the strength of L2 regularization, with smaller values implying weak regularization which allows the model to fit the training data more closely and risks overfitting, and larger values implying strong regularization which forces the model coefficients to shrink towards zero and reduces model complexity and overfitting. Intermediate values provide a spectrum of regularization strengths to balance between underfitting and overfitting. To find the right balance, we provide options with a range of magnitudes (0.0001, 0.05, 0.5, and 0.8). The choice of learning rate schedule affects the convergence behavior \citep{lecun2002efficient, bengio2012practical, bottou2018optimization}: constant learning schedule maintains a fixed step size for stable but potentially slow convergence; inverse scaling gradually reduces the step size to balance exploration and fine-tuning; and the adaptive learning schedule adjusts dynamically based on gradient changes to improve convergence efficiency. We use these three options in the set of hyperparameters. Table \ref{table:hyper} lists the hyperparameter values considered and those selected. We use the ``Adam" solver for model training, which is suitable for relatively large datasets like ours \citep{kingma2014adam}. We use mean squared error as the scoring metric, which has smooth differentiability, making it well suited for such gradient-based optimization.

\vspace{-4mm}

\begin{table}[hbt]
\scriptsize
\centering
\begin{tabular}{l c c c } \hline
\textbf{Grid Search Hyperparameter} & \textbf{Candidate Values} & \textbf{Selected Value (no RAS)} & \textbf{Selected Value (with RAS)}\\ \hline
 Hidden layer size & 1000, 2500,  4000 & 4000 & 4000\\ 
 Activation function & Identity, ReLU & ReLU & ReLU\\ 
 Regularization parameter & 0.0001, 0.05, 0.5, 0.8 & 0.05 & 0.0001\\ 
 Learning rate & Constant, Inverse Scaling, Adaptive & Constant & Constant\\ \hline
\end{tabular}
\caption{Candidate hyperparameter values for NN grid search.}
\label{table:hyper}
\end{table}

\subsection{Drone-Only Trajectory Model}\label{SA:DOSDTM} \label{sec:4.2}
In order to generate the training data for the NN predictor, we modify the lower-level model as follows. 

We remove the constraint \eqref{eq:sec3p2_oepr_time4} in order to use the total drone travel time as the objective function. In addition, we remove the truck-and-drone coordination constraints \eqref{eq:sec3p2_touchz1}-\eqref{eq:sec3p2_touchv1}. Finally, constraints \eqref{eq:sec3p2_flightd1}-\eqref{eq:sec3p2_flightd3} are replaced by constraint \eqref{eq:sec3p1_flightd1}, and constraint \eqref{eq:sec3p2_flightdv1}-\eqref{eq:sec3p2_flightdv3} are replaced by constraint \eqref{eq:sec3p1_flightdv1}.

\vspace{-8mm}

\begin{subequations}
\begin{align}
r_{i,t+1} = r_{it} + \Delta_f \cdot v_{it} + \frac{\Delta^2}{2} \cdot a_{it} & \qquad \forall i \in \{x,y\}, t \in \mathcal{T}^-_{f} \label{eq:sec3p1_flightd1}\\
v_{i,t+1} = v_{it} + \Delta_f \cdot a_{it} &\qquad \forall i \in \{x,y\}, t \in \mathcal{T}^-_{f} \label{eq:sec3p1_flightdv1}
\end{align}
\end{subequations}

\subsection{Three Solution Approaches to Estimate Drone Operation Time} \label{sec:4.3} 
We present three methods to calculate the drone travel time with varying degrees of simplicity.

\begin{enumerate}
\item Euclidean Distance-Based Method (K method): This method represents the state-of-the-art in drone travel time estimation for truck-and-drone planning. It computes the straight-line distances from the operation start node to the delivery node and the delivery node to the operation end node, and divides the sum of these two distances by an estimated drone speed (70 km/h) to estimate drone travel time.

\item Calibrated Euclidean Distance-Based Method (MK method): This method is an improved version of the state-of-the-art in drone time estimation for truck-and-drone planning. This method generates an offline training dataset using the process described in Section \ref{sec:4.1}, listing the coordinates of the start, delivery, and end nodes, for a training set of operations. We then run the drone-only model (from Section \ref{SA:DOSDTM}) on each operation to obtain the drone travel times. Next, we divide this travel time by the Euclidean distance-based drone travel time calculated by the K method to produce an operation-specific correction factor. We calculate the average value of this correction factor over all operations in the training set. For any new operation presented in the online phase, we calculate the drone travel time of the K method, and multiply it by this average correction factor from the training phase to get the final drone travel time estimate.

\item Neural Network Based Method (P method): This method trains an NN on the offline training set, and uses the trained NN to estimate drone travel time for any new operation presented in the online phase.
\end{enumerate}

\section{Computational Results}\label{section_computation}
We conducted computational experiments for two different scenarios. Scenario I approximates a truck-and-drone service area without tall buildings or mountains, that is, without restricted airspaces (RAS), and with uniformly distributed delivery nodes. Scenario II represents real-world settings with large regions where the truck and drone are not allowed (modeled as RAS), due to tall buildings or mountains, and with uniformly distributed delivery needs in the rest of the service area. The truck is assumed to travel along the shortest path between consecutive nodes avoiding the RAS, if any. In each scenario, we compare the duration and drone energy consumption of the entire truck-and-drone tour for the three solution approaches (K, MK, and P). We vary the number of delivery nodes (across 10, 20, 50, 75, 100, 175, 250) and the average speed of the truck (across 20, 30, 40, 50, 60, 70, 80 km/h). For each combination of scenario, number of nodes, and truck speed, we generate 100 instances with different sets of delivery locations and report average results. Each truck-and-drone tour plan in this section was computed in at most 189 seconds of runtime. 

\subsection{Scenario I: No Restricted Airspaces}\label{sec:computation_scenarioI}
We first quantify the benefits of drone use. Then we evaluate the benefits of our P method compared to benchmark approaches in terms of reduction in tour duration and, finally, in terms of the drone energy saved.

\begin{table}[hbt]
\scriptsize
\setlength{\tabcolsep}{2pt} 
\centering
\begin{tabular}{>{\centering\arraybackslash}p{2.8cm}| 
>{\centering\arraybackslash}p{0.65cm} 
>{\centering\arraybackslash}p{0.65cm} 
>{\centering\arraybackslash}p{0.65cm} 
>{\centering\arraybackslash}p{0.65cm} 
>{\centering\arraybackslash}p{0.65cm} 
>{\centering\arraybackslash}p{0.65cm} 
>{\centering\arraybackslash}p{0.65cm} || 
>{\centering\arraybackslash}p{0.65cm} 
>{\centering\arraybackslash}p{0.65cm} 
>{\centering\arraybackslash}p{0.65cm} 
>{\centering\arraybackslash}p{0.65cm} 
>{\centering\arraybackslash}p{0.65cm} 
>{\centering\arraybackslash}p{0.65cm} 
>{\centering\arraybackslash}p{0.65cm} }
\toprule
\multicolumn{1}{c|}{Metrics} 
& \multicolumn{7}{c||}{\textbf{Instances Where P Outperforms Truck-only TSP}} 
& \multicolumn{7}{c}{\textbf{Avg. \% Reduction in Duration by P method}} \\
\midrule
\diagbox[width=2.8cm]{\scriptsize Nodes}{\scriptsize Speed (km/h)} 
& 20 & 30 & 40 & 50 & 60 & 70 & 80 
& 20 & 30 & 40 & 50 & 60 & 70 & 80 \\
\midrule
10  & 100 & 99 & 100 & 100 & 99 & 99 & 99  
    & 25.28 & 23.49 & 21.22 & 17.67 & 14.00 & 12.41 & 11.74 \\
20  & 100 & 100 & 100 & 100 & 100 & 100 & 99  
    & 25.88 & 24.59 & 22.15 & 18.79 & 15.92 & 13.90 & 12.08 \\
50  & 100 & 100 & 100 & 100 & 100 & 98  & 100
    & 24.28 & 22.48 & 19.23 & 16.18 & 13.73 & 11.88 & 10.67 \\
75  & 100 & 100 & 100 & 100 & 100 & 100 & 100
    & 23.97 & 21.49 & 18.13 & 15.53 & 13.11 & 11.23 & 10.18 \\
100 & 100 & 100 & 100 & 100 & 100 & 100 & 100
    & 24.12 & 21.02 & 17.90 & 15.13 & 12.86 & 10.90 & 9.75 \\
175 & 100 & 100 & 100 & 100 & 100 & 100 & 100
    & 22.65 & 19.46 & 16.18 & 13.86 & 11.83 & 9.86 & 8.83 \\
250 & 100 & 100 & 100 & 100 & 100 & 100 & 100
    & 21.98 & 18.48 & 15.29 & 13.02 & 11.09 & 9.04 & 8.35 \\
\bottomrule
\end{tabular}
\vspace{0.2cm}
\caption{Comparing the P method with the Truck-only TSP for Scenario I (without RAS).  
Columns 2-8 list the number of instances (out of 100) where the P method strictly outperforms the Truck-only TSP in terms of the tour duration.  
Columns 9-15 list the average percentage reduction in duration by the P method compared to Truck-only TSP.}
\label{table:Computational_studies_figure1}
\end{table}

\begin{table}[hbt]
\scriptsize
\setlength{\tabcolsep}{2pt} 
\centering
\begin{tabular}{>{\centering\arraybackslash}p{2.8cm}| 
>{\centering\arraybackslash}p{0.65cm} 
>{\centering\arraybackslash}p{0.65cm} 
>{\centering\arraybackslash}p{0.65cm} 
>{\centering\arraybackslash}p{0.65cm} 
>{\centering\arraybackslash}p{0.65cm} 
>{\centering\arraybackslash}p{0.65cm} 
>{\centering\arraybackslash}p{0.65cm} || 
>{\centering\arraybackslash}p{0.65cm} 
>{\centering\arraybackslash}p{0.65cm} 
>{\centering\arraybackslash}p{0.65cm} 
>{\centering\arraybackslash}p{0.65cm} 
>{\centering\arraybackslash}p{0.65cm} 
>{\centering\arraybackslash}p{0.65cm} 
>{\centering\arraybackslash}p{0.65cm} }
\toprule
\multicolumn{1}{c|}{Metrics} 
& \multicolumn{7}{c||}{\textbf{Avg. \% Reduction in Duration by P vs. K}} 
& \multicolumn{7}{c}{\textbf{Avg. \% Reduction in Duration by P vs. MK}} \\
\midrule
\diagbox[width=2.8cm]{\scriptsize Nodes}{\scriptsize Speed (km/h)} 
& 20 & 30 & 40 & 50 & 60 & 70 & 80
& 20 & 30 & 40 & 50 & 60 & 70 & 80 \\ 
\midrule
10  & 0.15 & 0.73 & 2.30 & 4.32 & 6.78 & 5.66 & 3.97
    & 0.02 & 0.26 & 0.21 & 0.95 & 0.38 & 0.41 & 0.69 \\
20  & 0.25 & 1.60 & 4.22 & 10.56 & 11.16 & 8.89 & 7.09
    & 0.08 & 0.50 & 1.37 & 2.05 & 1.28 & 0.97 & 1.21 \\
50  & 0.52 & 3.84 & 11.19 & 18.56 & 19.78 & 16.06 & 13.67
    & 0.30 & 2.25 & 6.72 & 7.21 & 3.12 & 3.19 & 2.70 \\
75  & 0.83 & 5.54 & 14.67 & 22.89 & 24.00 & 19.70 & 16.07
    & 0.47 & 4.23 & 10.34 & 10.36 & 5.61 & 5.93 & 4.23 \\
100 & 1.41 & 7.51 & 18.51 & 26.06 & 27.43 & 21.86 & 18.57
    & 1.13 & 5.97 & 13.24 & 14.04 & 7.21 & 7.05 & 6.46 \\
175 & 3.73 & 13.96 & 26.23 & 33.75 & 34.08 & 27.47 & 23.80
    & 3.33 & 12.23 & 21.04 & 20.65 & 12.40 & 11.69 & 10.55 \\
250 & 6.82 & 20.64 & 33.23 & 39.80 & 38.77 & 31.21 & 28.20
    & 6.47 & 18.74 & 28.19 & 26.01 & 16.17 & 15.42 & 14.45 \\
\bottomrule
\end{tabular}
\vspace{0.2cm}
\caption{Comparing the P method with the K and MK methods on total duration for Scenario I (without RAS). Columns 2-8 (respectively, columns 9-15) list the average percentage reduction by the P method compared to the K (respectively, MK) method. Each row shows a different number of nodes and each column a different truck speed.}
\label{table:Computational_studies_figure2}
\end{table}

\begin{table}[hbt]
\scriptsize
\setlength{\tabcolsep}{2pt} 
\centering
\begin{tabular}{>{\centering\arraybackslash}p{2.8cm}| 
>{\centering\arraybackslash}p{0.75cm} 
>{\centering\arraybackslash}p{0.75cm} 
>{\centering\arraybackslash}p{0.75cm} 
>{\centering\arraybackslash}p{0.75cm} 
>{\centering\arraybackslash}p{0.75cm} 
>{\centering\arraybackslash}p{0.75cm} 
>{\centering\arraybackslash}p{0.75cm} || 
>{\centering\arraybackslash}p{0.75cm} 
>{\centering\arraybackslash}p{0.75cm} 
>{\centering\arraybackslash}p{0.75cm} 
>{\centering\arraybackslash}p{0.75cm} 
>{\centering\arraybackslash}p{0.75cm} 
>{\centering\arraybackslash}p{0.75cm} 
>{\centering\arraybackslash}p{0.75cm} }
\toprule
\multicolumn{1}{c|}{Metrics} 
& \multicolumn{7}{c||}{\textbf{Avg. Reduction (\%) in DEC by P vs. K}} 
& \multicolumn{7}{c}{\textbf{Avg. Reduction (\%) in DEC by P vs. MK}} \\
\midrule
\diagbox[width=2.8cm]{\scriptsize Nodes}{\scriptsize Speed (km/h)}
& 20 & 30 & 40 & 50 & 60 & 70 & 80
& 20 & 30 & 40 & 50 & 60 & 70 & 80 \\
\midrule
10  & 1.82 & 3.82 & 9.38 & 13.98 & 22.30 & 22.65 & 13.26
    & -0.80 & -1.52 & 3.25 & -1.70 & -14.66 & -4.24 & -6.43 \\
20  & 2.91 & 8.17 & 17.49 & 25.24 & 24.60 & 21.00 & 17.91
    & 2.16 & 5.80 & 12.96 & 10.34 & -1.00 & -3.55 & -3.92 \\
50  & 8.20 & 21.16 & 32.36 & 36.78 & 35.74 & 31.13 & 25.05
    & 7.71 & 19.70 & 27.73 & 24.62 & 12.26 & 8.98 & 6.84 \\
75  & 12.99 & 27.88 & 37.68 & 42.86 & 40.92 & 33.90 & 30.31
    & 12.47 & 26.63 & 33.43 & 31.05 & 19.11 & 14.34 & 14.16 \\
100 & 18.01 & 33.20 & 42.01 & 46.23 & 44.17 & 37.24 & 33.36
    & 17.79 & 31.87 & 37.78 & 35.58 & 22.34 & 18.10 & 17.53 \\
175 & 27.87 & 43.51 & 51.48 & 54.58 & 52.33 & 44.53 & 41.42
    & 27.50 & 42.21 & 47.64 & 43.91 & 32.86 & 28.33 & 26.93 \\
250 & 36.13 & 50.86 & 57.81 & 60.40 & 57.78 & 49.60 & 45.03
    & 35.80 & 49.54 & 54.14 & 49.82 & 38.84 & 34.31 & 31.76 \\
\bottomrule
\end{tabular}
\vspace{0.2cm}
\caption{Comparing the P method with the K and MK methods on DEC for Scenario I (without RAS). Columns 2-8 (respectively, columns 9-15) list the average percentage reduction by the P method compared to the K (respectively, MK) method. Each row shows a different number of nodes and each column a different truck speed.}
\label{table:Computational_studies_figure3}
\end{table}

\subsubsection{Benefits of Using a Drone:}\label{sec:comp_benefitofdrone}
In the conventional delivery logistics, all packages are delivered by the truck on a TSP tour that visits all nodes. To evaluate the extent to which drone usage can reduce tour duration, we compare the durations of truck-only and truck-and-drone tours. We solve the problem of minimizing the truck-only tour duration using the TSP heuristics used in the first step of our overall solution approach, while we minimize the truck-and-drone tour duration using the aforementioned P method. For a given number of nodes and truck speed, let the total durations of the truck-and-drone and truck-only tours for the $i^\text{th}$ instance (of the 100 instances) be $\tau_i^P$ and $\tau _i^T$, respectively. We calculate how frequently the truck-and-drone solution outperforms the truck-only solution as $\sum_{i \in \{1,\cdots, 100\}} I (\tau _i^T > \tau _i^P)$ and report it in columns 2-8 in Table \ref{table:Computational_studies_figure1} and the left sub-figure of Figure \ref{fig:s1_pvt2} in Appendix \ref{section_appendix_comp_plot}. We also calculate the average reduction in total duration achieved by the truck-and-drone tour compared to the truck-only tour as $\frac{1}{100}\sum_{i \in \{1,\cdots, 100\}} \frac{\tau _i^T- \tau _i^P}{ \tau _i^T }$ and report it in columns 9-15 in Table \ref{table:Computational_studies_figure1} and the right sub-figure of Figure \ref{fig:s1_pvt2} in Appendix \ref{section_appendix_comp_plot}.

From these results, we observe that in 43 of the 49 combinations of the number of nodes and truck speed, the truck-and-drone solution has strictly less duration than the truck-only solution for \textit{all} 100 instances. Across all 4,900 instances, the truck-and-drone solution has strictly less duration than the truck-only solution in 4,893 (that is, 99.86\%) instances. The truck-only solution performs as well as the truck-and-drone P method only in a handful of instances with a small number of nodes (50 or less), especially at high truck speeds (70 km/h or higher). High truck speeds dilute the speed advantage of the drone, making the drone less useful. Columns 9-15 in Table \ref{table:Computational_studies_figure1} show that the drone reduces the tour duration, on average, by 16.38\%.

Columns 2-8 in Table \ref{table:Computational_studies_figure1} show that the duration reduction occurs more frequently at lower truck speeds (which amplify the drone's speed advantage) and more nodes. On the other hand, columns 9-15 in Table \ref{table:Computational_studies_figure1} show that the greatest average improvement also occurs at lower truck speeds (as expected), but at \textit{fewer} nodes. At first glance, this observation is somewhat surprising but can be explained as follows. As the number of nodes increases, there are more opportunities for the drone to improve the truck-only TSP solution because there are more candidate drone nodes, and hence many more truck-and-drone operations to explore. So, there is a higher likelihood of finding at least some node that, when assigned to a drone, can reduce the overall travel time. Thus, there is a greater chance to improve the tour using a drone.
Therefore, more nodes leads to more frequent duration reduction. But with fewer nodes, the denominator in the calculation of percentage reduction (i.e., the truck-only tour duration) decreases, making any improvement appear larger as a percentage. So, in case of fewer nodes, the delegation of a subset of delivery nodes to the drone yields a higher percentage reduction in tour duration, although the improvement frequency itself is lower.

Overall, the truck-and-drone tour optimized with our P method consistently and substantially outperforms the truck-only tour, especially at low to medium truck speeds, which are common in large congested cities.

\subsubsection{Benefits of Our Solution Approach:}\label{sec:comp_benefitofsa}
Table \ref{table:Computational_studies_figure2} and Figure \ref{fig:s1_pvkvmk_time} (in Appendix \ref{section_appendix_comp_plot}) show the average reduction in truck-and-drone tour duration achieved by the P method compared to the K and MK methods. As before, we solve 100 instances for each combination of the number of nodes and truck speed. For the $i^\text{th}$ instance, let $\tau_i^P$, $\tau _i^K$ and $\tau_i^{MK}$ be the truck-and-drone tour durations obtained by the P, K and MK method, respectively. Then the average reduction in truck-and-drone tour duration achieved by the P method compared to the K and MK methods is calculated as $\frac{1}{100}\sum_{i \in \{1,\cdots, 100\}} \frac{\tau _i^K- \tau _i^P}{ \tau _i^K }$ and $\frac{1}{100}\sum_{i \in \{1,\cdots, 100\}} \frac{\tau _i^{MK}- \tau _i^P}{ \tau _i^{MK} }$.

We observe that in each of the 49 combinations, the P method outperforms both the K and MK methods in terms of tour duration, on average by 15.06\% and 7.13\%, respectively. The duration reduction increases monotonically as the number of nodes increases. This higher percentage of duration reduction for a larger number of nodes can be intuitively explained as follows. As demonstrated later (in Section \ref{sec:networkstructure}), K and MK methods tend to overestimate the drone's duration reduction ability and thus overuse the drone. So, with more nodes, the total time wasted by one vehicle waiting for the other accumulates and amounts to an especially longer tour duration for the K and MK methods, relative to the P method. The P method does not suffer from this problem given its accurate estimate of the drone travel time. As the number of nodes increases, the number of operations grows faster than linearly, and the total tour duration grows slower than linearly (a well-known result from the TSP literature), together resulting in a higher percentage reduction for the P method compared to the K and MK methods. In practice, this observation implies that as the demand for e-commerce delivery increases, a better estimate of drone travel time brings even greater benefits.

In terms of truck speed, the duration reduction peaks at around 60 km/h compared to the K method, and around 50 km/h compared to the MK method. This suggests that our P method adds the most value at intermediate truck speeds, which can be explained by recognizing the following two opposite factors. On the one hand, at high truck speeds, the usefulness of a drone diminishes regardless of the chosen optimization method, resulting in smaller differences between methods. On the other hand, at low truck speeds, since the drone travels much faster than the truck, the different methods result in a similar set of operations being chosen, lowering the relative benefits of the P method: the large differential between truck and drone speeds implies that the operations favored by the K and MK methods under their overestimation of the drone's ability are also favored under the more realistic drone speed estimates by the P method. The net effect of these two opposite factors amounts to the highest benefits of the P method at intermediate truck speeds. 

Overall, the P method consistently outperforms the K and MK methods, with the most impressive benefits in the more practical parameter ranges, namely, at more delivery locations and intermediate truck speeds.

\subsubsection{Drone Energy Consumption Reduction:}\label{sec:comp_droneenergy}
As mentioned in Section \ref{sec:TaDTPM}, we use drone energy consumption (DEC) minimization as a secondary objective to break the ties if multiple tours have the same duration. Let $\gamma_i^P$, $\gamma_i^K$, and $\gamma_i^{MK}$, respectively, be the DEC of the solutions obtained for the $i^\text{th}$ instance using the P, K and MK methods. We calculate the average drone energy savings as $\frac{1}{100}\sum_{i \in \{1,\cdots, 100\}} \frac{\gamma _i^K- \gamma _i^P}{ \tau _i^K }$ and $\frac{1}{100}\sum_{i \in \{1,\cdots, 100\}} \frac{\gamma _i^{MK}- \gamma _i^P}{ \gamma _i^{MK} }$, and report them in Table \ref{table:Computational_studies_figure3} and Figure \ref{fig:s1_pvkvmk_energy} (in Appendix \ref{section_appendix_comp_plot}).

The P method reduces the average DEC for all 49 combinations, compared to the K method, and for 40 of the 49 combinations compared to the MK method. On average, the P method produces a 31.61\%  and 20.01\% lower DEC, compared to the K and MK methods, respectively. Similar to the pattern observed for the duration reduction, DEC savings also increase with more nodes, as the inefficient assignment of deliveries to the two vehicles causes more hovering time for the drone waiting for the truck, wasting more drone energy. As expected, under such suboptimal assignments, DEC accumulates faster than linearly with the increasing number of deliveries. Compared to the reduction in tour duration, the DEC savings peak earlier: at 40 km/h compared to K and 30 km/h compared to the MK method, because higher truck speeds undermine the DEC savings more directly than the tour duration savings. Overall, as with the duration reduction, our P method produces substantial DEC reduction, especially with more nodes and low to medium truck speeds.

\subsection{Scenario II: With Restricted Airspaces}

\begin{table}[hbt]
\scriptsize
\setlength{\tabcolsep}{2pt} 
\centering
\begin{tabular}{>{\centering\arraybackslash}p{2.8cm}| 
>{\centering\arraybackslash}p{0.7cm} 
>{\centering\arraybackslash}p{0.7cm} 
>{\centering\arraybackslash}p{0.7cm} 
>{\centering\arraybackslash}p{0.7cm} 
>{\centering\arraybackslash}p{0.7cm} 
>{\centering\arraybackslash}p{0.7cm} 
>{\centering\arraybackslash}p{0.7cm} || 
>{\centering\arraybackslash}p{0.7cm} 
>{\centering\arraybackslash}p{0.7cm} 
>{\centering\arraybackslash}p{0.7cm} 
>{\centering\arraybackslash}p{0.7cm} 
>{\centering\arraybackslash}p{0.7cm} 
>{\centering\arraybackslash}p{0.7cm} 
>{\centering\arraybackslash}p{0.7cm} }
\toprule
\multicolumn{1}{c|}{Metrics} 
& \multicolumn{7}{c||}{\textbf{Instances Where P Outperforms Truck-only TSP}} 
& \multicolumn{7}{c}{\textbf{Avg. Reduction (\%) in Duration by P method}} \\
\midrule
\diagbox[width=2.8cm]{\scriptsize Nodes}{\scriptsize Speed (km/h)} 
& 20 & 30 & 40 & 50 & 60 & 70 & 80
& 20 & 30 & 40 & 50 & 60 & 70 & 80 \\
\midrule
10  & 88  & 85  & 81  & 73  & 67  & 61  & 57 
    & 14.19 & 13.98 & 11.06 & 7.46 & 4.02 & 2.64 & 1.50 \\
20  & 98  & 98  & 96  & 91  & 85  & 82  & 83 
    & 18.57 & 17.28 & 14.94 & 11.55 & 8.81 & 7.43 & 6.22 \\
50  & 100 & 100 & 100 & 98  & 95  & 90  & 93 
    & 18.18 & 16.43 & 14.07 & 10.93 & 8.45 & 6.81 & 6.07 \\
75  & 100 & 100 & 100 & 100 & 98  & 98  & 96 
    & 19.83 & 17.42 & 14.81 & 11.52 & 8.80 & 7.30 & 6.36 \\
100 & 100 & 100 & 100 & 99  & 100 & 98  & 99 
    & 19.79 & 17.64 & 14.46 & 11.11 & 8.81 & 6.71 & 5.59 \\
175 & 100 & 100 & 98  & 100 & 100 & 98  & 96 
    & 18.57 & 15.61 & 13.07 & 10.34 & 8.30 & 6.72 & 5.22 \\
250 & 100 & 100 & 100 & 97  & 97  & 98  & 93 
    & 17.54 & 14.49 & 11.57 & 8.76 & 6.89 & 5.85 & 4.06 \\
\bottomrule
\end{tabular}
\vspace{0.2cm}
\caption{Comparing the P method with the Truck-only TSP for Scenario II (with RAS).  
Columns 2-8 list the number of instances (out of 100) where the P method strictly outperforms the Truck-only TSP in terms of the tour duration.  
Columns 9-15 list the average percentage reduction in duration by the P method compared to Truck-only TSP.}
\label{table:Computational_studies_figure4}
\end{table}

\begin{table}[hbt]
\scriptsize
\setlength{\tabcolsep}{2pt} 
\centering
\begin{tabular}{>{\centering\arraybackslash}p{2.8cm}| 
>{\centering\arraybackslash}p{0.7cm} 
>{\centering\arraybackslash}p{0.7cm} 
>{\centering\arraybackslash}p{0.7cm} 
>{\centering\arraybackslash}p{0.7cm} 
>{\centering\arraybackslash}p{0.7cm} 
>{\centering\arraybackslash}p{0.7cm} 
>{\centering\arraybackslash}p{0.7cm} || 
>{\centering\arraybackslash}p{0.7cm} 
>{\centering\arraybackslash}p{0.7cm} 
>{\centering\arraybackslash}p{0.7cm} 
>{\centering\arraybackslash}p{0.7cm} 
>{\centering\arraybackslash}p{0.7cm} 
>{\centering\arraybackslash}p{0.7cm} 
>{\centering\arraybackslash}p{0.7cm} }
\toprule
\multicolumn{1}{c|}{Metrics} 
& \multicolumn{7}{c||}{\textbf{Avg. Reduction (\%) in Duration by P vs. K}} 
& \multicolumn{7}{c}{\textbf{Avg. Reduction (\%) in Duration by P vs. MK}} \\
\midrule
\diagbox[width=2.8cm]{\scriptsize Nodes}{\scriptsize Speed (km/h)}
& 20 & 30 & 40 & 50 & 60 & 70 & 80
& 20 & 30 & 40 & 50 & 60 & 70 & 80 \\ 
\midrule
10  & 0.27 & 0.87 & 1.67 & 4.45 & 6.82 & 7.24 & 5.05
    & 0.02 & -0.05 & -0.47 & 0.27 & -0.43 & -0.18 & 0.05 \\
20  & 0.45 & 1.74 & 5.04 & 8.74 & 12.88 & 9.79 & 6.93
    & 0.20 & 0.68 & 1.20 & 1.69 & -0.03 & 0.07 & -0.49 \\
50  & 1.39 & 4.87 & 11.63 & 18.87 & 21.34 & 16.47 & 12.09
    & 0.92 & 3.25 & 7.54 & 7.40 & 2.43 & 1.11 & 1.23 \\
75  & 2.34 & 7.44 & 16.46 & 24.11 & 25.39 & 18.73 & 15.34
    & 1.91 & 6.05 & 12.07 & 11.45 & 3.61 & 2.64 & 2.65 \\
100 & 3.81 & 11.87 & 21.40 & 28.39 & 29.34 & 21.66 & 17.20
    & 3.38 & 9.63 & 16.65 & 15.02 & 6.16 & 4.23 & 4.33 \\
175 & 8.19 & 20.11 & 31.89 & 37.50 & 37.58 & 28.85 & 24.33
    & 7.72 & 18.30 & 26.42 & 24.82 & 12.66 & 10.68 & 9.63 \\
250 & 12.89 & 27.81 & 38.89 & 44.32 & 43.01 & 33.52 & 29.13
    & 12.39 & 25.97 & 33.85 & 30.52 & 17.69 & 16.44 & 14.02 \\
\bottomrule
\end{tabular}
\vspace{0.2cm}
\caption{Comparing the P method with the K and MK methods on total duration for Scenario II (with RAS). Columns 2-8 (respectively, columns 9-15) list the average percentage reduction by the P method compared to the K (respectively, MK) method. Each row shows a different number of nodes and each column a different truck speed.}
\label{table:Computational_studies_figure5}
\end{table}

\begin{table}[hbt]
\scriptsize
\setlength{\tabcolsep}{2pt} 
\centering
\begin{tabular}{
  >{\centering\arraybackslash}p{2.8cm}|
  >{\centering\arraybackslash}p{0.7cm}
  >{\centering\arraybackslash}p{0.7cm}
  >{\centering\arraybackslash}p{0.7cm}
  >{\centering\arraybackslash}p{0.7cm}
  >{\centering\arraybackslash}p{0.7cm}
  >{\centering\arraybackslash}p{0.7cm}
  >{\centering\arraybackslash}p{0.7cm}
  ||
  >{\centering\arraybackslash}p{0.7cm}
  >{\centering\arraybackslash}p{0.7cm}
  >{\centering\arraybackslash}p{0.7cm}
  >{\centering\arraybackslash}p{0.7cm}
  >{\centering\arraybackslash}p{0.7cm}
  >{\centering\arraybackslash}p{0.7cm}
  >{\centering\arraybackslash}p{0.7cm}
}
\toprule
\multicolumn{1}{c|}{Metrics} 
& \multicolumn{7}{c||}{\textbf{Average reduction (\%) in DEC by P vs. K}} 
& \multicolumn{7}{c}{\textbf{Average reduction (\%) in DEC by P vs. MK}} \\
\midrule
\diagbox[width=2.8cm]{\scriptsize Nodes}{\scriptsize Speed (km/h)}
& 20 & 30 & 40 & 50 & 60 & 70 & 80
& 20 & 30 & 40 & 50 & 60 & 70 & 80 \\
\midrule
10  & 1.32 & 2.80 & 6.10 & 14.40 & 21.35 & 21.47 & 11.10
    & -1.15 & -0.23 & -0.37 & 3.19 & -7.81 & -2.70 & -1.38 \\
20  & 1.80 & 4.77 & 11.98 & 19.15 & 25.62 & 22.44 & 15.99
    & 1.07 & 2.93 & 6.69 & 8.59 & -0.12 & -3.00 & -3.81 \\
50  & 6.93 & 16.49 & 25.27 & 31.94 & 33.35 & 28.23 & 20.10
    & 6.39 & 14.65 & 20.74 & 19.26 & 5.73 & 1.54 & 1.63 \\
75  & 11.37 & 22.70 & 31.78 & 37.33 & 36.59 & 30.65 & 23.12
    & 10.72 & 21.00 & 27.59 & 24.78 & 9.00 & 5.96 & 4.70 \\
100 & 15.92 & 29.33 & 38.59 & 41.70 & 41.34 & 32.23 & 25.24
    & 15.32 & 27.21 & 33.77 & 29.02 & 14.27 & 7.46 & 6.80 \\
175 & 26.98 & 40.88 & 49.32 & 51.80 & 49.58 & 40.35 & 33.57
    & 26.46 & 39.36 & 44.64 & 39.72 & 22.46 & 17.53 & 14.85 \\
250 & 34.85 & 49.09 & 55.87 & 58.55 & 55.89 & 46.12 & 39.65
    & 34.46 & 47.60 & 51.55 & 46.81 & 30.46 & 26.22 & 22.07 \\
\bottomrule
\end{tabular}
\vspace{0.2cm}
\caption{Comparing the P method with the K and MK methods on DEC for Scenario II (with RAS). Columns 2-8 (respectively, columns 9-15) list the average percentage reduction by the P method compared to the K (respectively, MK) method. Each row shows a different number of nodes and each column a different truck speed.}
\label{table:Computational_studies_figure6}
\end{table}

\subsubsection{Benefits of Using a Drone:}
Table \ref{table:Computational_studies_figure4} and Figure \ref{fig:s2_pvt2} (in Appendix \ref{section_appendix_comp_plot}) demonstrate the benefits of using a drone in Scenario II. Across 4,900 instances, 4,586 (that is, 93.59\%) result in a strict improvement in tour duration. Similar to Scenario I, columns 2-8 in Table \ref{table:Computational_studies_figure4} show that the duration reduction occurs more frequently at lower truck speeds and more nodes, and columns 9-15 in Table \ref{table:Computational_studies_figure4} show that the average percentage reduction is also higher at lower truck speeds. However, unlike Scenario I (where columns 9-15 in Table \ref{table:Computational_studies_figure1} show a higher average percentage reduction at fewer nodes), columns 9-15 in Table \ref{table:Computational_studies_figure4} show that, with more nodes, the benefit of using the drone first increases and then decreases. This is because at a lower number of nodes, additional nodes present more chances for the drone to help with the delivery (a much stronger effect than in Scenario I, as shown in columns 2-8 in Table \ref{table:Computational_studies_figure4} vs. Table \ref{table:Computational_studies_figure1}), which contributes to improving the average benefit of the drone. However, as the number of nodes grows further, the denominator in the percentage reduction becomes larger, resulting in a smaller percentage reduction in duration.

\subsubsection{Benefits of Our Solution Approach:}
Table \ref{table:Computational_studies_figure5} and Figure \ref{fig:s2_pvkvmk_time} (in Appendix \ref{section_appendix_comp_plot}) quantify the average percentage reduction in tour duration achieved by the P method compared to the K and MK methods, in Scenario II. The P method outperforms the K method in all 49 combinations and outperforms the MK method in 43 of the 49 combinations. The average percentage reduction in duration is 16.74\% with respect to the K method and 7.90\% relative to the MK method, which align closely (and are slightly better) than Scenario I. As with Scenario I, the duration reduction tends to increase as the number of nodes grows, but peaks at intermediate truck speeds of around 50-60 km/h compared to the K method and at around 40-50 km/h compared to the MK method. Overall, we find that, even with the incorporation of RAS, the duration reduction benefits are comparable to those without RAS and exhibit similar broad patterns.

\subsubsection{Drone Energy Consumption Reduction:}
Table \ref{table:Computational_studies_figure6} and Figure \ref{fig:s2_pvkvmk_energy} (in Appendix \ref {section_appendix_comp_plot}) present the DEC reduction using the P method compared to the K and MK methods in Scenario II. The average reduction in DEC is 28.43\% and 15.79\% compared to the K and MK methods—numbers that are similar, albeit slightly lower, than in Scenario I. As the node count increases, the advantage of the P method generally increases. In terms of truck speed, DEC savings also peak around 50-60 km/h compared to the K method and around 40-50 km/h compared to the MK method. In summary, the qualitative and quantitative findings related to the benefits of P method remain largely the same regardless of whether the service area has RAS. 

\section{Practical Case Studies and Managerial Insights} \label{sec6}
To demonstrate the effectiveness of our mathematical and computational approach in real-world settings, we conducted case studies in four major US city centers. We directly use the NN predictor trained in Scenario I in Section \ref{section_computation} to run the experiments in this section, demonstrating its robustness across different settings.

We obtained road networks and building locations from the OpenStreetMap NetworkX Python package (OSMNX) \citep{boeing2024modeling}. For each selected region in the city (5 km $\times$ 5 km, centered on downtown), we identified all buildings and set the probability of having a delivery node at a building to be proportional to the volume of each building. From real-world data in these regions, we observed that the total footprint of all \textit{tall} buildings (taller than the drone's maximum allowed flying altitude) constituted less than 1\% of the total area (5 km $\times$ 5 km). Thus, in these case studies, there are effectively no restricted airspaces. We sampled 200 delivery locations for a typical day, repeating 50 times to ensure statistical robustness and produce average results for each city. We assumed a maximum drone speed of 70 km/h \citep{faa_part107} and an average truck speed of 40 km/h (i.e., 25 mi/h), common in busy downtowns. Each truck-and-drone tour plan in our practical case studies was computed in at most 38 seconds of runtime.

For convenient visualization, Section \ref{sec:networkstructure} uses a smaller (50-node) example from Section \ref{section_computation} to demonstrate the key differences between the truck-and-drone tours generated using our P method and those using the two algorithmic benchmarks (K and MK methods). Then Section \ref{sec:benefitsofoptimization} quantifies the degree to which the P method outperforms the K and MK methods, as well as the truck-only TSP solutions for the real-world case studies described at the beginning of this section. Finally, in Section \ref{sec:driversofimprovement}, we identify the main attributes of the solutions generated by our approach that drive these improvements for real-world case studies.

\subsection{Tour Structure} \label{sec:networkstructure}
\begin{figure}[ht]
    \centering
    \includegraphics[scale=0.55]{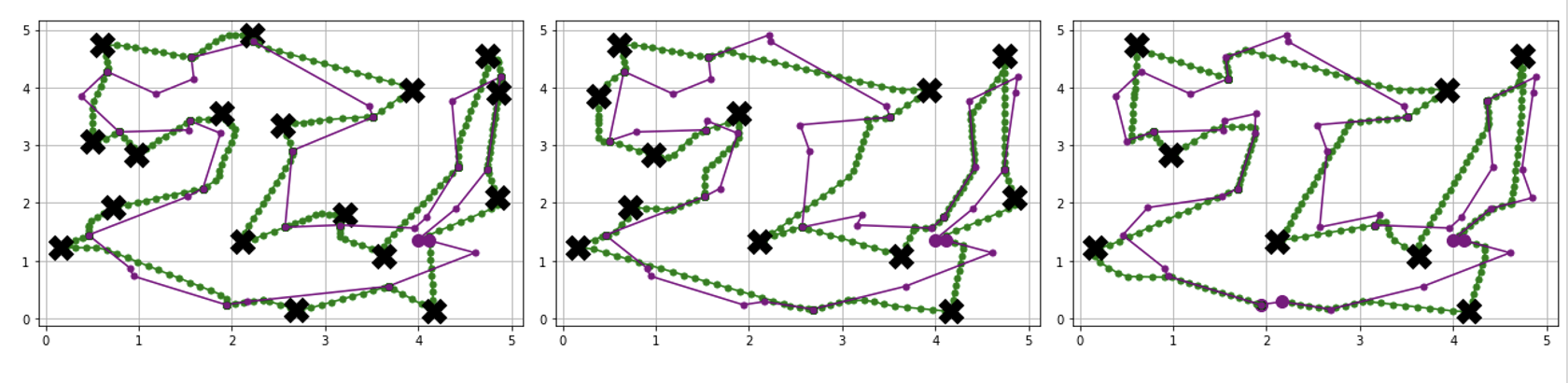}
    \caption{Example tour solutions using K (left), MK (middle) and P (right) methods with the same instance.}
    \label{fig:ex_tour_3_method}
\end{figure}

We now analyze common structures in our truck-and-drone solutions compared to the benchmarks, to understand why we see relative improvements, and to identify the key drivers of improvement explained in Section \ref{sec:driversofimprovement}. Figure \ref{fig:ex_tour_3_method} shows the optimal tours obtained using the K, MK, and P methods for a representative example instance with 50 uniformly generated delivery locations without RAS. Black crosses represent drone delivery locations, green dotted lines represent drone trajectories, purple dots represent truck delivery locations, and purple lines represent truck paths. Across the three solutions, the drone typically serves delivery nodes located near the truck path, on either side, but requires inconvenient detours to the truck path to serve these nodes via truck. In such cases, drone deliveries are more efficient by avoiding such detours. 

\vspace{-5mm}

\vspace{-5mm}

We observe that in the K and MK solutions, which assume Euclidean and modified Euclidean drone travel, the drone's effectiveness in reducing the tour duration is overestimated, resulting in more drone nodes (17 and 12 in this example), compared to the solution from the P method (which has eight nodes here). Moreover, the drone nodes selected in the P solution are often subsets of the drone nodes in the K and MK solutions. These patterns are broadly seen beyond this illustrative example, which motivates Section \ref{sec:driversofimprovement}.

\subsection{Benefits of Optimization} \label{sec:benefitsofoptimization}
\begin{table}[htbp]
\centering
\scriptsize
\caption{Case studies using straight-line distances for truck travel (T stands for the Truck-only TSP tour).}
\label{table:case_euclidean}
\begin{tabular}{llcccccc}
\toprule
\textbf{Metric} & \textbf{Stat} & \textbf{Boston} & \textbf{Chicago} & \textbf{Philadelphia} & \textbf{Los Angeles} & \textbf{Average} \\
\midrule
\multirow{3}{*}{K vs. T duration} 
    & Mean & -89.24\% & -82.75\% & -37.56\% & -29.05\% & -59.65\% \\
    & 95\% CI & {[}-92.02\%, -86.47\%{]} & {[}-85.18\%, -80.32\%{]} & {[}-39.47\%, -35.65\%{]} & {[}-30.93\%, -27.17\%{]} & {[}-61.90\%, -57.40\%{]} \\
    & Count   & 0 & 0 & 0 & 0 & 0 \\
\midrule
\multirow{3}{*}{MK vs. T duration}
    & Mean & -72.42\% & -63.29\% & -27.25\% & -20.21\% & -45.79\% \\
    & 95\% CI & {[}-75.34\%, -69.51\%{]} & {[}-66.06\%, -60.51\%{]} & {[}-29.09\%, -25.41\%{]} & {[}-22.11\%, -18.31\%{]} & {[}-48.15\%, -43.43\%{]} \\
    & Count   & 0 & 0 & 0 & 0 & 0 \\
\midrule
\multirow{3}{*}{P vs. T duration}
    & Mean & 9.87\% & 10.64\% & 7.53\% & 7.67\% & 8.93\% \\
    & 95\% CI & {[}8.13\%, 11.61\%{]} & {[}8.89\%, 12.40\%{]} & {[}6.28\%, 8.77\%{]} & {[}6.38\%, 8.96\%{]} & {[}7.42\%, 10.43\%{]} \\
    & Count & 47 & 46 & 48 & 46 & 46.75 \\
\midrule
\multirow{3}{*}{P vs. K duration}
    & Mean & 52.31\% & 51.02\% & 32.64\% & 28.28\% & 41.06\% \\
    & 95\% CI & {[}51.42\%, 53.21\%{]} & {[}49.93\%, 52.11\%{]} & {[}31.49\%, 33.80\%{]} & {[}26.92\%, 29.64\%{]} & {[}39.94\%, 42.19\%{]} \\
    & Count   & 50 & 50 & 50 & 50 & 50 \\
\midrule
\multirow{3}{*}{P vs. MK duration}
    & Mean & 47.64\% & 45.16\% & 27.19\% & 22.92\% & 35.73\% \\
    & 95\% CI & {[}46.69\%, 48.60\%{]} & {[}44.02\%, 46.29\%{]} & {[}26.00\%, 28.38\%{]} & {[}21.18\%, 24.66\%{]} & {[}34.47\%, 36.98\%{]} \\
    & Count   & 50 & 50 & 50 & 50 & 50 \\
\midrule
\multirow{3}{*}{P vs. K DEC}
    & Mean & 71.56\% & 70.37\% & 60.05\% & 56.19\% & 64.54\% \\
    & 95\% CI & {[}71.09\%, 72.03\%{]} & {[}69.86\%, 70.89\%{]} & {[}59.49\%, 60.60\%{]} & {[}55.59\%, 56.79\%{]} & {[}64.01\%, 65.08\%{]} \\
    & Count   & 50 & 50 & 50 & 50 & 50 \\
\midrule
\multirow{3}{*}{P vs. MK DEC}
    & Mean & 68.49\% & 66.58\% & 55.90\% & 52.04\% & 60.75\% \\
    & 95\% CI & {[}67.81\%, 69.17\%{]} & {[}65.86\%, 67.30\%{]} & {[}55.21\%, 56.58\%{]} & {[}51.34\%, 52.74\%{]} & {[}60.05\%, 61.45\%{]} \\
    & Count   & 50 & 50 & 50 & 50 & 50 \\
\bottomrule
\end{tabular}
\end{table}

\begin{table}[htbp]
\centering
\scriptsize
\caption{Case studies using road network-based distances for truck travel (T stands for the Truck-only TSP tour).}
\label{table:case_roadnetwork}
\begin{tabular}{llcccccc}
\toprule
\textbf{Metric} & \textbf{Stat} & \textbf{Boston} & \textbf{Chicago} & \textbf{Philadelphia} & \textbf{Los Angeles} & \textbf{Average} \\
\midrule
\multirow{3}{*}{K vs. T duration} 
    & Mean    & 26.46\% & 24.44\% & 18.58\% & 14.12\% & 20.90\% \\
    & 95\% CI & {[}25.56\%, 27.37\%{]} & {[}23.47\%, 25.42\%{]} & {[}17.70\%, 19.47\%{]} & {[}13.00\%, 15.25\%{]} & {[}19.93\%, 21.88\%{]} \\
    & Count   & 50 & 50 & 50 & 50 & 50 \\
\midrule
\multirow{3}{*}{MK vs. T duration}
    & Mean    & 27.12\% & 25.13\% & 19.02\% & 14.65\% & 21.48\% \\
    & 95\% CI & {[}26.17\%, 28.07\%{]} & {[}24.20\%, 26.05\%{]} & {[}18.13\%, 19.92\%{]} & {[}13.57\%, 15.72\%{]} & {[}20.52\%, 22.44\%{]} \\
    & Count   & 50 & 50 & 50 & 50 & 50 \\
\midrule
\multirow{3}{*}{P vs. T duration}
    & Mean    & 28.48\% & 26.25\% & 20.27\% & 17.26\% & 23.06\% \\
    & 95\% CI & {[}27.63\%, 29.33\%{]} & {[}25.39\%, 27.12\%{]} & {[}19.39\%, 21.14\%{]} & {[}16.38\%, 18.13\%{]} & {[}22.20\%, 23.93\%{]} \\
    & Count   & 50 & 50 & 50 & 50 & 50 \\
\midrule
\multirow{3}{*}{P vs. K duration}
    & Mean    & 2.70\% & 2.31\% & 2.02\% & 3.52\% & 2.64\% \\
    & 95\% CI & {[}1.92\%, 3.48\%{]} & {[}1.31\%, 3.32\%{]} & {[}1.19\%, 2.85\%{]} & {[}2.33\%, 4.71\%{]} & {[}1.69\%, 3.59\%{]} \\
    & Count   & 39 & 38 & 40 & 41 & 39.5 \\
\midrule
\multirow{3}{*}{P vs. MK duration}
    & Mean    & 1.81\% & 1.43\% & 1.49\% & 2.94\% & 1.92\% \\
    & 95\% CI & {[}1.10\%, 2.53\%{]} & {[}0.45\%, 2.42\%{]} & {[}0.71\%, 2.28\%{]} & {[}1.77\%, 4.11\%{]} & {[}1.00\%, 2.84\%{]} \\
    & Count   & 38 & 33 & 40 & 38 & 37.25 \\
\midrule
\multirow{3}{*}{P vs. K DEC}
    & Mean    & 14.02\% & 17.88\% & 13.92\% & 21.20\% & 16.76\% \\
    & 95\% CI & {[}13.16\%, 14.88\%{]} & {[}16.88\%, 18.88\%{]} & {[}13.06\%, 14.78\%{]} & {[}19.89\%, 22.51\%{]} & {[}15.75\%, 17.76\%{]} \\
    & Count   & 50 & 50 & 50 & 50 & 50 \\
\midrule
\multirow{3}{*}{P vs. MK DEC}
    & Mean    & 12.75\% & 16.64\% & 13.16\% & 20.60\% & 15.78\% \\
    & 95\% CI & {[}11.80\%, 13.69\%{]} & {[}15.70\%, 17.58\%{]} & {[}12.29\%, 14.02\%{]} & {[}19.28\%, 21.92\%{]} & {[}14.77\%, 16.80\%{]} \\
    & Count   & 50 & 50 & 50 & 50 & 50 \\
\bottomrule
\end{tabular}
\end{table}

\begin{figure}[hbt]
    \centering
    \includegraphics[scale=0.5]{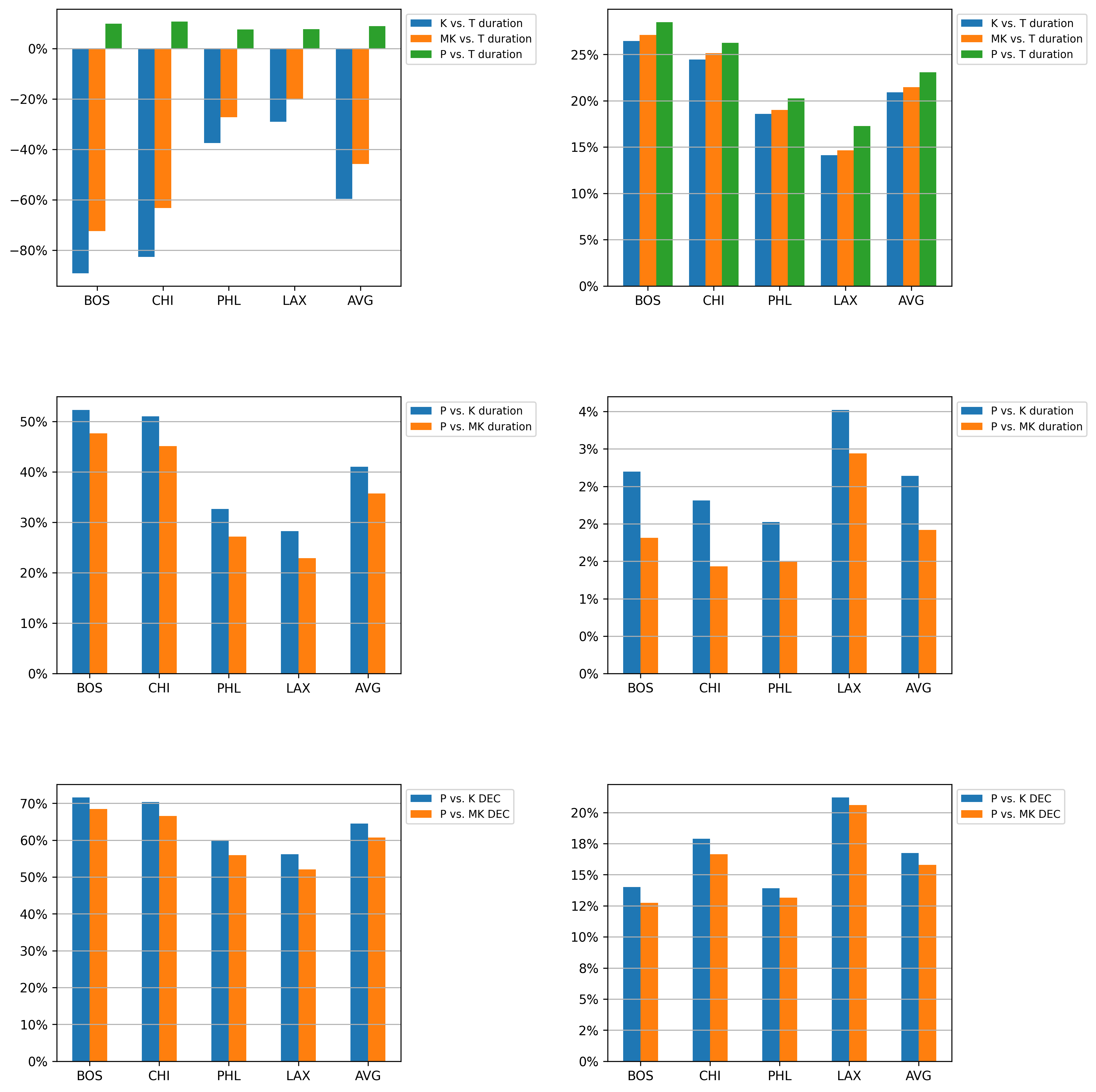}
    \caption{Case study results assuming straight-line (left) and road network-based (right) truck travels.}
    \label{fig:case_compare_plot}
\end{figure}

In this section, we quantify the degree of improvement achieved by the P, K, and MK methods, compared to the truck-only TSP solution, and also the improvement achieved by the P method compared to the K and MK methods. We consider two different assumptions about truck travel: straight-line Euclidean travel (an assumption used in most of the existing literature, ensuring easier comparison) and travel along the actual road network (a more realistic assumption, ensuring more practicality). Tables \ref{table:case_euclidean} and \ref{table:case_roadnetwork} and Figure \ref{fig:case_compare_plot} summarize the improvements in terms of tour duration (in the first five main rows of the tables) and DEC (in the last two main rows). Each main row presents the mean and the 95\% confidence interval of the percentage improvement, as well as the count of how many of the 50 instances show a positive improvement.

Across four cities and both truck travel assumptions, our P method always performs the best in terms of both tour duration and DEC. Depending on the specific structures of the demand distributions and road networks in each downtown area, the actual duration reduction and the DEC savings vary. Across four cities and two truck travel assumptions, the P method reduces tour duration by 2.02\%-52.31\% compared to the K method and by 1.43\%-47.64\% compared to the MK method. In addition, the P method reduces DEC by 13.92\%-71.56\% compared to the K method and by 12.75\%-68.49\% compared to the MK method. P method also outperforms the T method (that is, the truck-only TSP tour) by 7.53\%-28.48\% in terms of the tour duration, across the four cities and two truck travel assumptions. Under road network-based truck travel, the K and MK methods also outperform truck-only tour by 14.12\%-27.12\% demonstrating that, in this case, the drone adds value regardless of the optimization method used in truck-and-drone tour planning. However, as shown in Table \ref{table:case_euclidean}, under the straight-line Euclidean truck travel assumption, the tour duration actually worsens significantly for the K and MK methods compared to the T method. This highlights the fact that, in some cases, inaccurate estimates of drone travel times can completely nullify the benefits of a drone.

Although Table \ref{table:case_roadnetwork} shows a smaller improvement achieved by the P method compared to the K, MK and T methods than in Table \ref{table:case_euclidean}, the 95\% confidence intervals demonstrate that all improvements are statistically significant. Moreover, these improvements are also practically very significant, from the point of view of the e-commerce companies. Amazon delivers more than 10 million customer packages per day \citep{contimod2025amazon}, and each driver makes 180 stops and delivers 250-300 packages a day \citep{capitalone2024amazonlogistics}, which corresponds to more than 33,333 truck drivers in service each day. Assuming a total cost of operation for a truck of \$90.78/hour \citep{fisher2023cost}, the resulting 1.92\% average reduction in duration using the P method, even compared to the better of the two benchmarks (the MK method) and even under the assumption of road network-based truck travel, translates into daily savings of more than \$200,000.

In summary, we find that our P method consistently and substantially outperforms all benchmarks under all tested assumptions, providing significant benefits both from a statistical and a practical standpoint.

\subsection{Drivers of Improvement} \label{sec:driversofimprovement}
We further analyze the results of case studies in which the truck travels along the road network in the downtown areas of Boston (BOS), Philadelphia (PHL), Chicago (CHI), and Los Angeles (LAX), to identify three main drivers of the improvement achieved by our P method compared to the benchmarks. These drivers include more truck delivery nodes, judicious reconfiguration of drone operations to retain most high-quality drone nodes, and large drone travel savings at the expense of small increases in truck travel distances.

\begin{figure}[hbt]
    \centering
    \begin{subfigure}[t]{0.45\textwidth} 
        \centering
        \includegraphics[width=\textwidth]{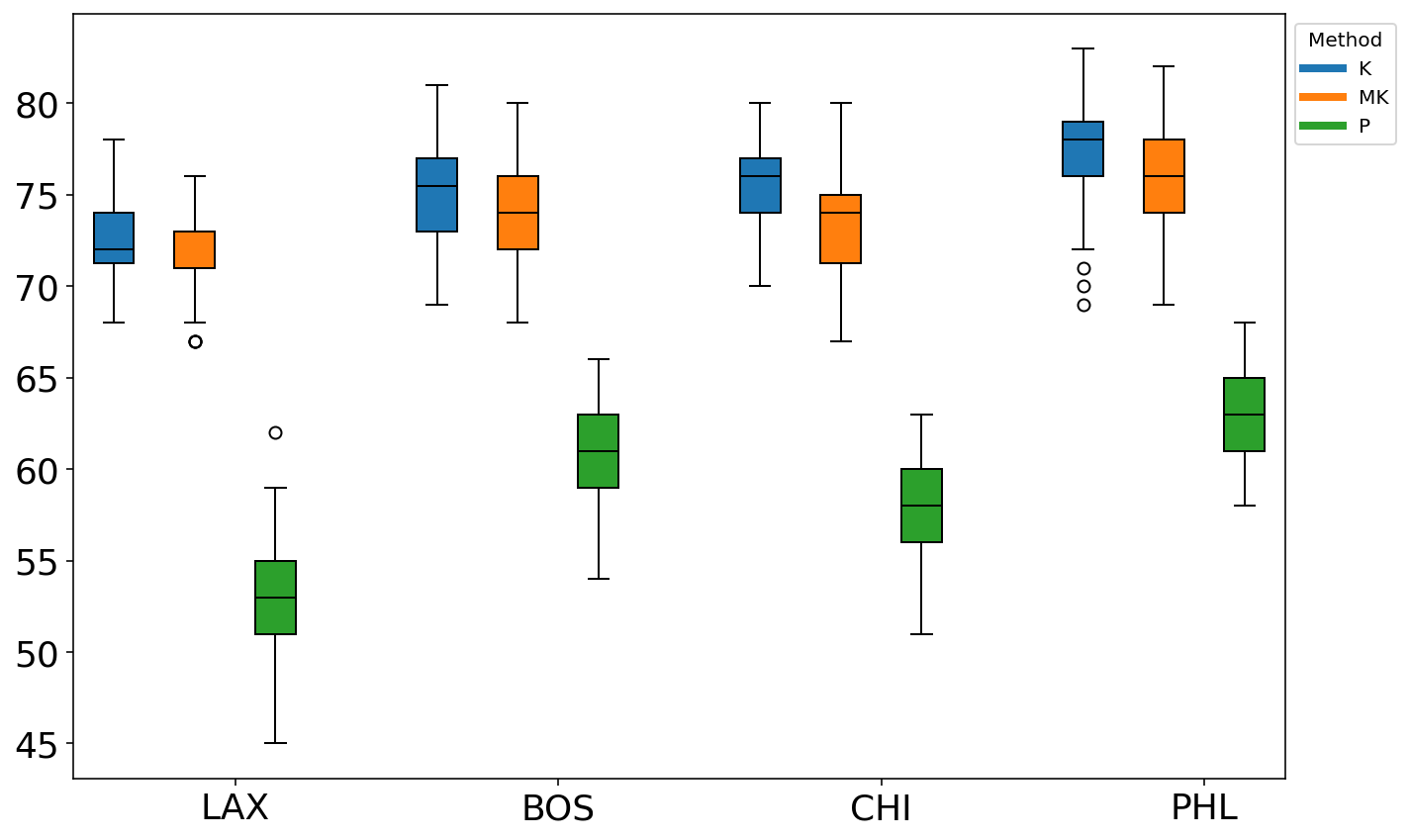}
        \caption{Number of drone nodes.}
        \label{fig:num_drone_nodes}
    \end{subfigure}
    \hspace{0.02\textwidth} 
    \begin{subfigure}[t]{0.45\textwidth} 
        \centering
        \includegraphics[width=\textwidth]{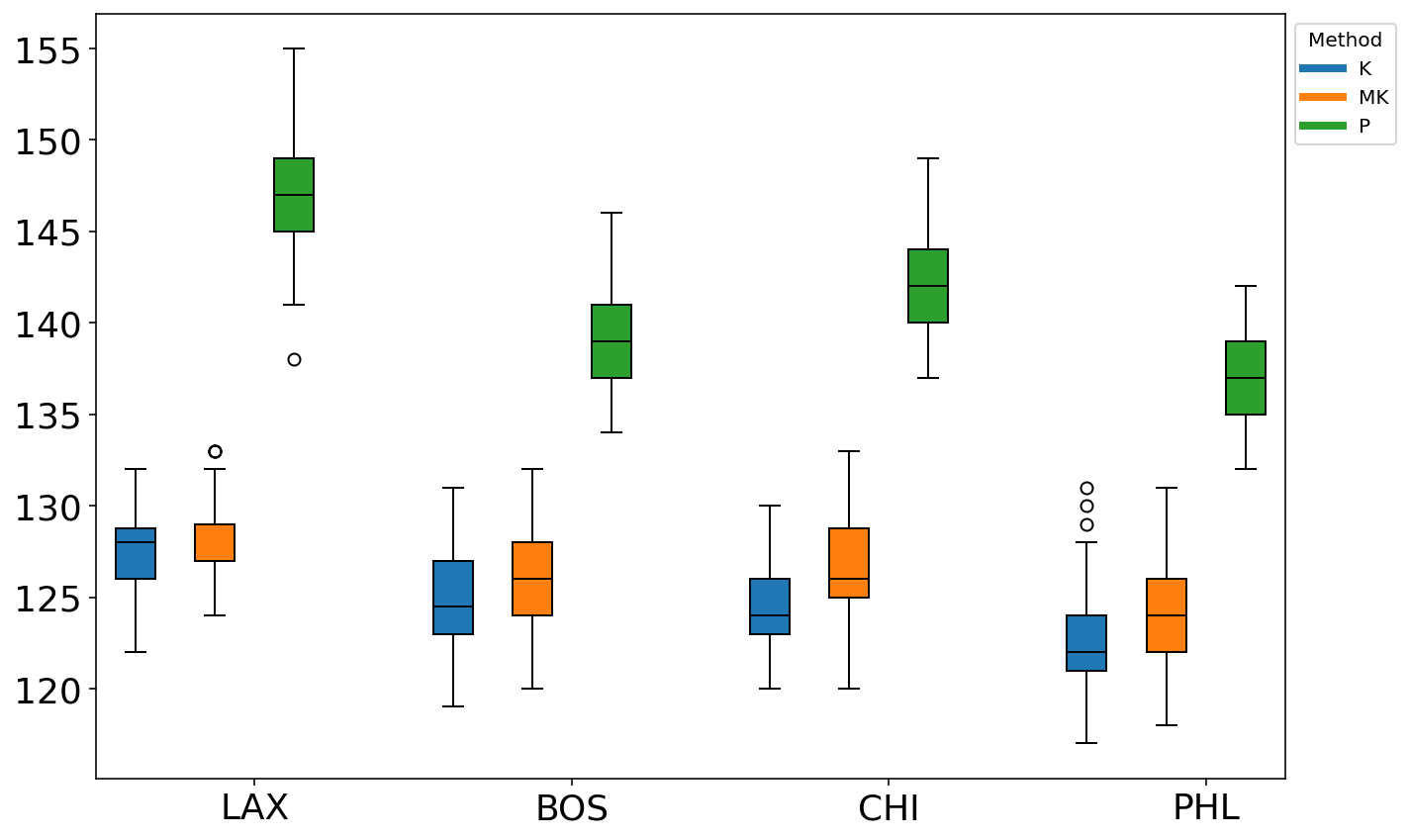}
    \caption{Number of truck nodes.}
    \label{fig:num_truck_nodes}
    \end{subfigure}
    \caption{Number of drone nodes (left) and truck nodes (right) averaged across the 50 instances.}
    \label{fig:comparison_v2}
\end{figure}

\textbf{More truck delivery nodes}: The K method, used extensively by previous studies in truck-and-drone tour planning, assumes Euclidean travel distances for the drone, while the MK method, developed by us, partially alleviates this shortcoming using a multiplier to adjust the drone travel time. Because of their reliance on Euclidean distances, we found that both K and MK methods tend to overestimate the drone's ability for fast deliveries. They both tend to use more drone operations than the P method, end up making the truck wait long periods for the drone to return, and thus increase the tour duration. Figure \ref{fig:comparison_v2} quantifies this effect by plotting the number of nodes served by the drone (in the left sub-figure) and the truck (in the right sub-figure) for each of the four cities and the three methods, namely K, MK and P. Figure \ref{fig:comparison_v2} shows that, across the four cities, the P method reduces the number of drone nodes by 21.7\% and 20.3\%, resulting in 13.0\% and 11.8\% more truck delivery nodes for the P method, compared to the K and MK methods, respectively. This overuse of drone nodes degrades the performance of the K and MK methods, and avoiding this pitfall in turn is an important driver of the superior performance of the P method.

\begin{figure}[hbt]
    \centering
    \begin{subfigure}[t]{0.45\textwidth} 
        \centering
        \includegraphics[width=\textwidth]{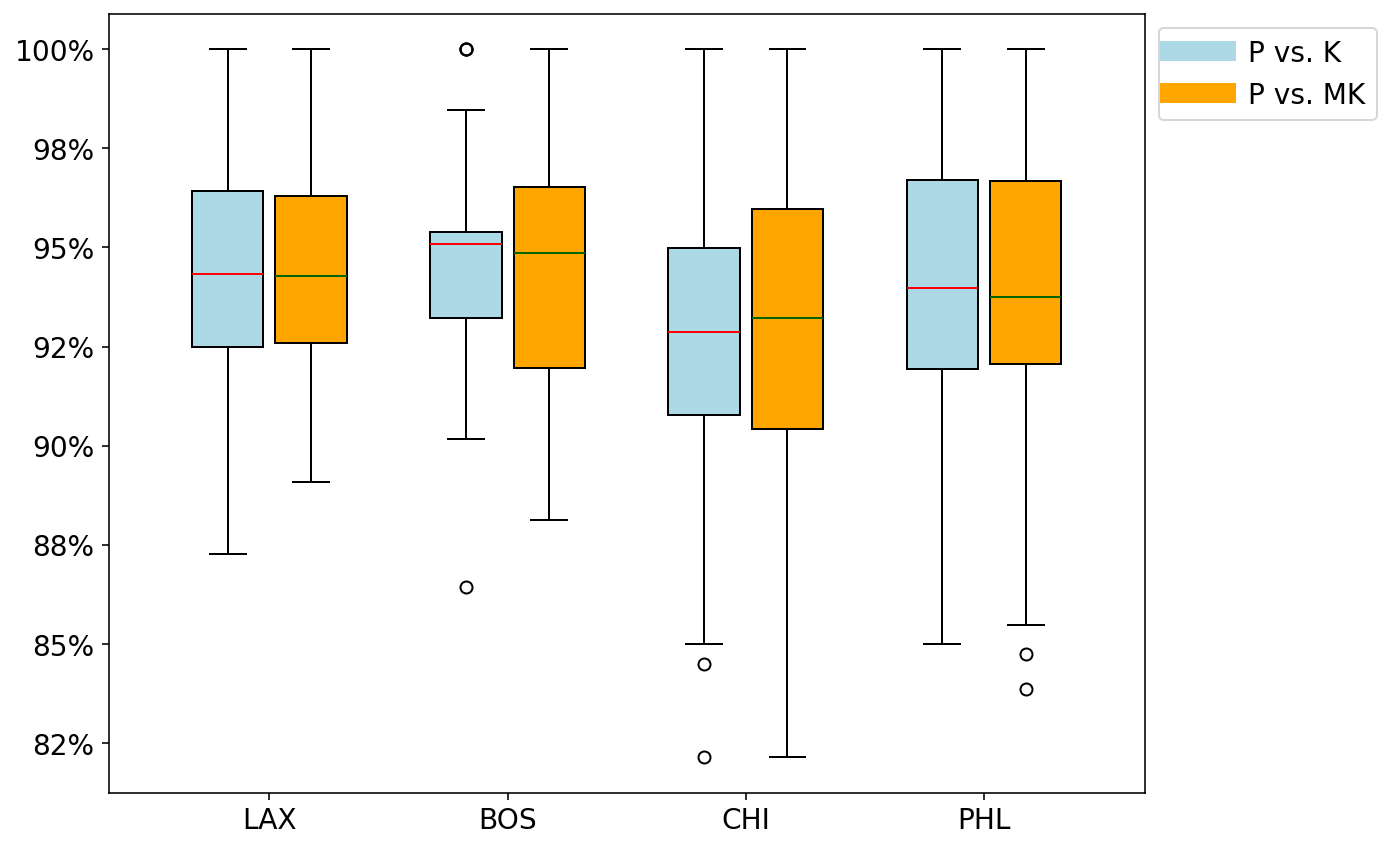}
        \caption{Percentage of P method's drone nodes that are drone nodes in the benchmark method.}
        \label{fig:common_drone_nodes}
    \end{subfigure}
    \hspace{0.02\textwidth} 
    \begin{subfigure}[t]{0.45\textwidth} 
        \centering
        \includegraphics[width=\textwidth]{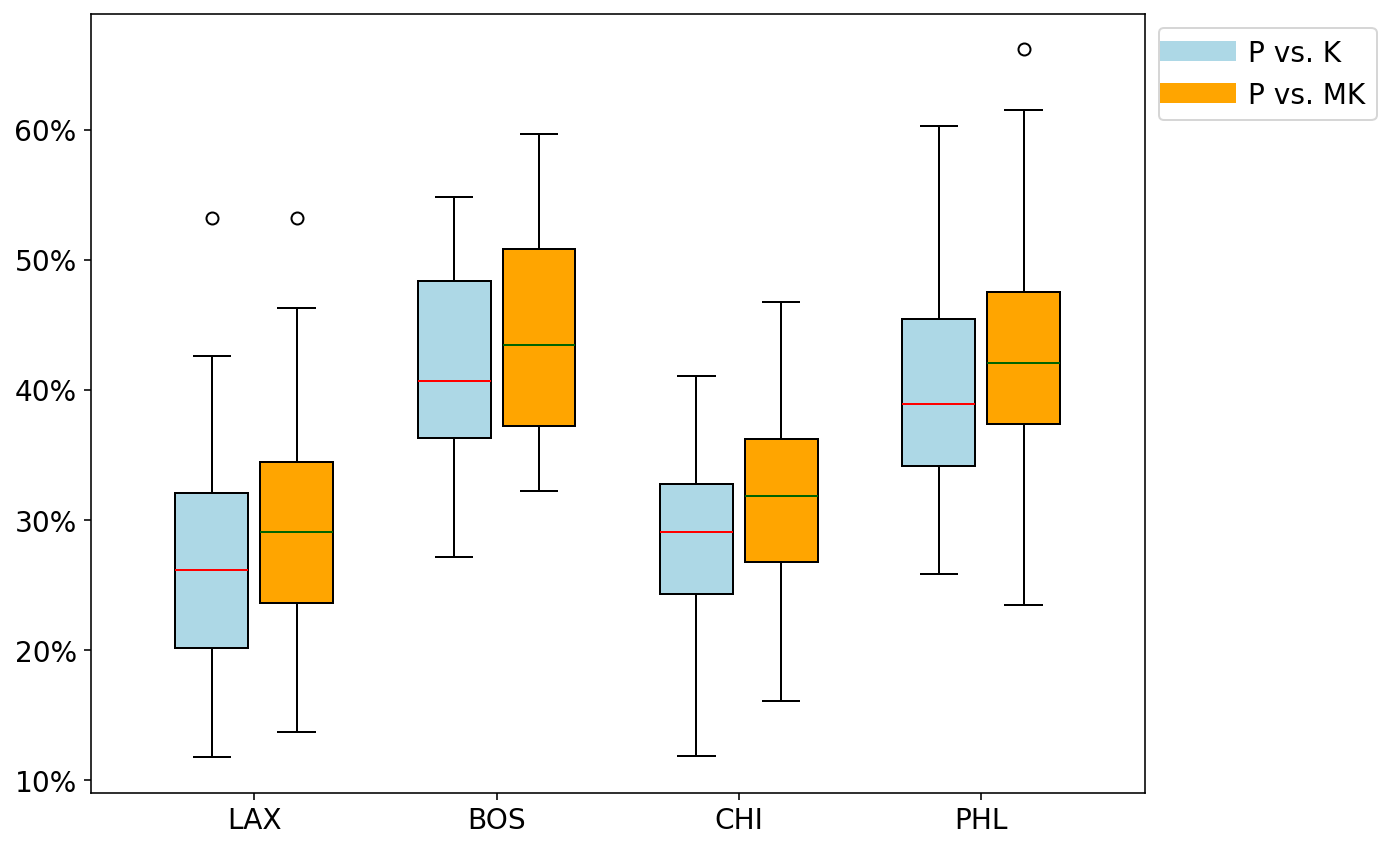}
    \caption{Percentage of P method's drone operations that are drone operations in the benchmark method.}
    \label{fig:common_drone_oper_wotruck}
    \end{subfigure}
    \caption{Percentage of P method's drone nodes and drone operations that overlap with those in benchmarks.}
    \label{fig:comparison_v3}
\end{figure}

\textbf{Judicious reconfiguration of drone operations to retain high-quality drone nodes}: In addition to avoiding overreliance on drone operations, the P method succeeds in retaining almost all high-quality drone nodes. Figure \ref{fig:common_drone_nodes} shows the number of drone nodes in the P method's solution that are also assigned to the drone in a benchmark method (K and MK) as a percentage of the total number of drone nodes in the P method's solution. The right sub-figure analogously shows the same metrics, but for drone operations rather than drone nodes, where a drone operation is defined as a 3-tuple consisting of the operation start node, drone node, and operation end node. Across all 200 instances in the four cities, 94\% of the drone nodes in the P method's solution are also found to be drone nodes in each of the K and MK methods' solutions.

However, a very different picture emerges in terms of overlaps in drone operations across different methods, as seen in Figure \ref{fig:common_drone_oper_wotruck}. Across all instances in the four cities, only 34.3\% and 37.2\% of the drone operations in the P method's solution are in the solutions of the K and MK methods, respectively. The contrast is even more striking if we focus on matching the entire operation, i.e., matching operation start node, operation end node, drone node and the entire sequence of truck-only nodes in an operation. Across the 200 instances in these four cities, \textit{none} of the operations chosen by the P method completely matches any of those chosen by either of the benchmarks (K or MK). Overall, across the 200 instances, although almost all ($\sim 94\%$) of the drone nodes selected by the P method are those in the benchmarks, only a third (34.3\%-37.2\%) of the drone operations match those in the benchmark solutions, and none of the operations match entirely. This suggests that the P method judiciously reconfigures the drone operations by changing the operation start/end nodes and truck visit sequence to generate more efficient operations, without relinquishing the high-quality drone nodes, which is a major driver of improvement obtained by the P method.

\begin{figure}[hbt]
    \centering
    \begin{subfigure}[t]{0.45\textwidth} 
        \centering
        \includegraphics[width=\textwidth]{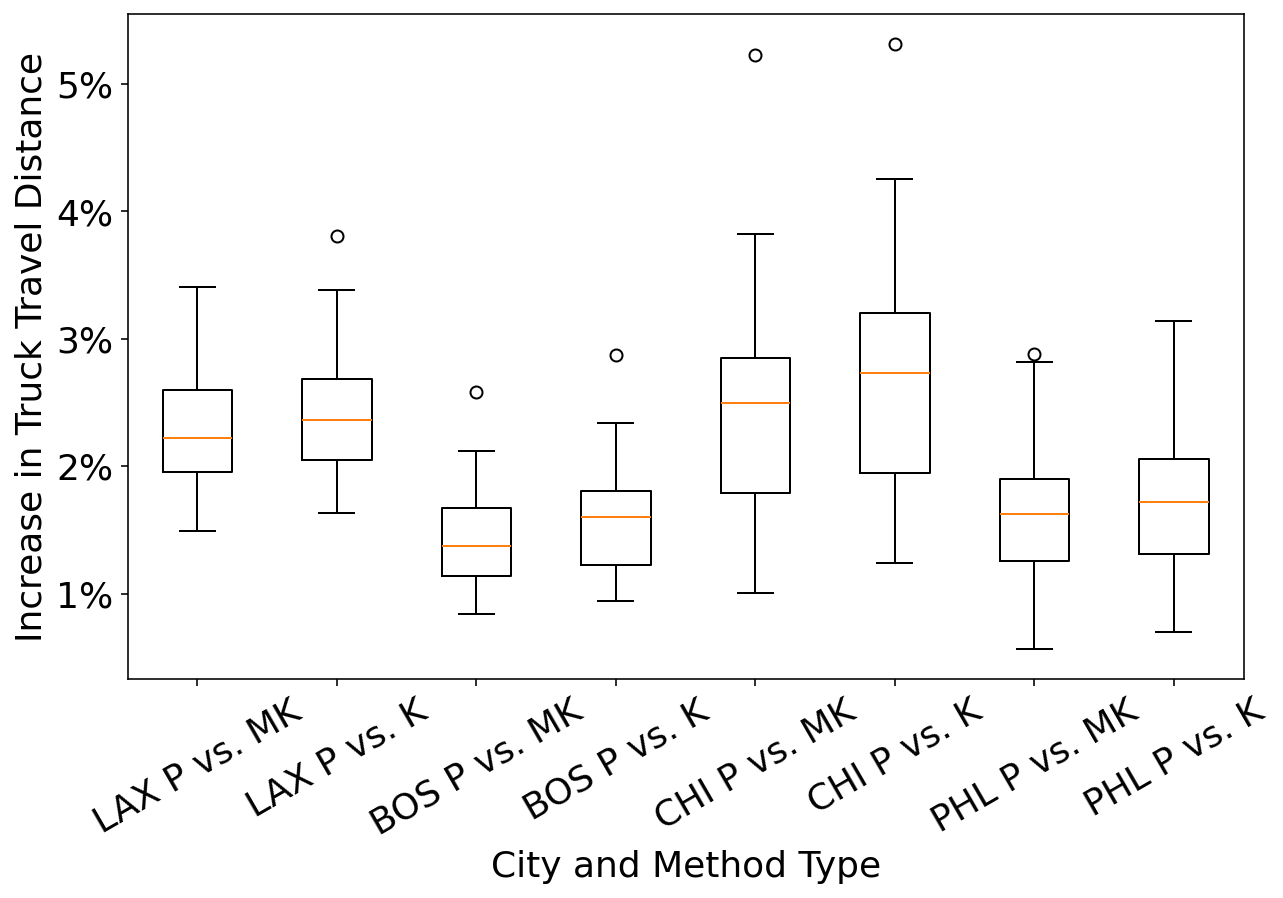}
        \caption{Percentage increase in truck travel distances.}
        \label{fig:truck_dist_increase}
    \end{subfigure}
    \hspace{0.02\textwidth} 
    \begin{subfigure}[t]{0.45\textwidth} 
        \centering
        \includegraphics[width=\textwidth]{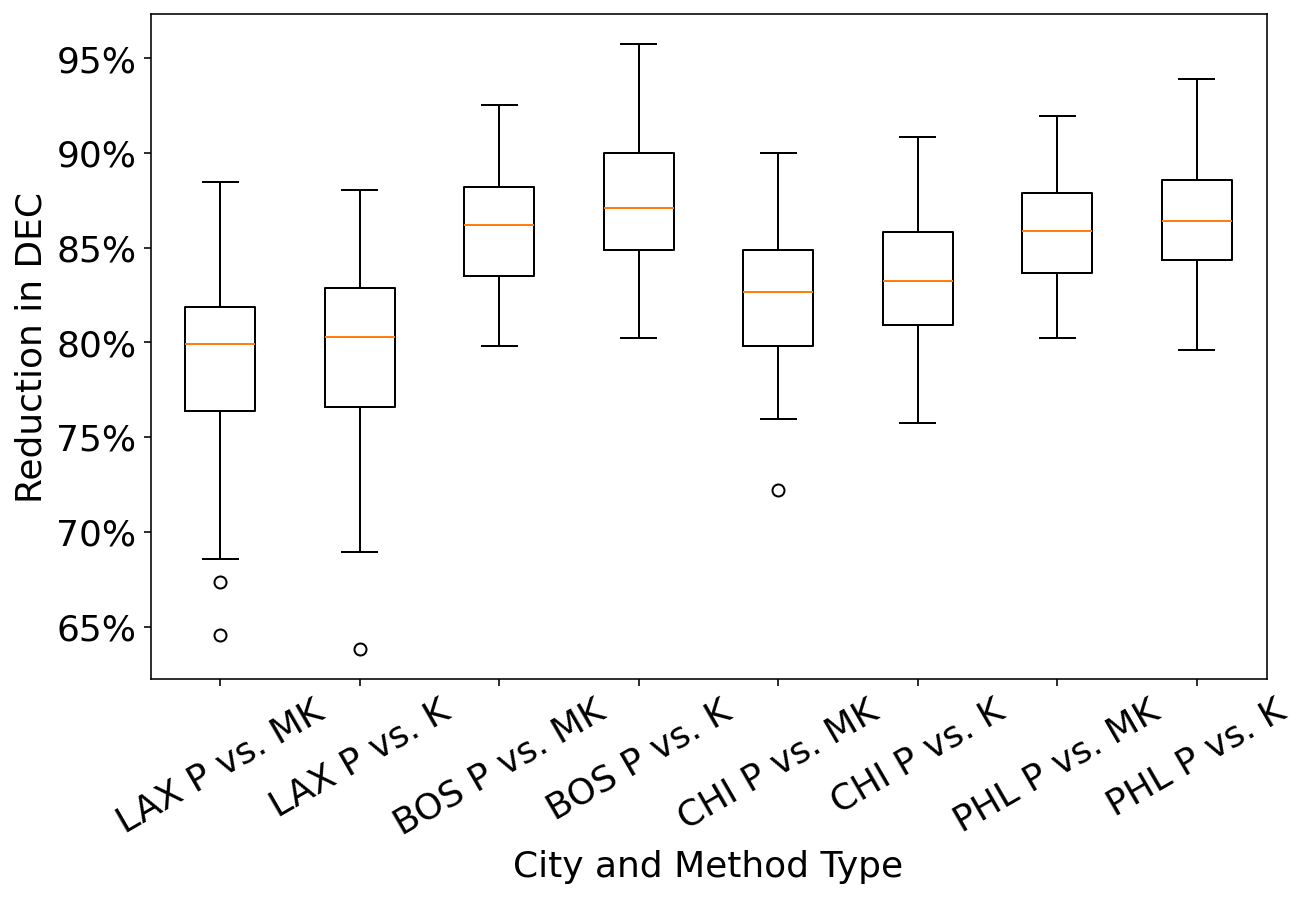}
        \caption{Percentage reduction in drone energy cost.}
        \label{fig:drone_energy_reduction}
    \end{subfigure}
    \caption{Increase in truck travel distance and DEC savings achieved by the P method relative to the benchmarks.}
    \label{fig:comparison_v1}
\end{figure}

\textbf{Large drone travel savings at the expense of small increases in truck travel distances}: A third key driver of improvement is the ability of the P method's solutions to achieve large reductions in drone travel in exchange of slight increases in truck travel. Figure \ref{fig:comparison_v1} shows that, compared to the K and MK methods, the P method increases the total distance traveled by the truck by only 2.1\% and 1.9\% of the the truck-only TSP length, but significantly reduces the drone energy consumption by 16.8\% and 15.8\%, respectively.

\section{Conclusion}
We formulate a nonlinear, multi-objective, truck-and-drone model that integrates drone flight physics within the combinatorial complexity of the truck routing problem. We propose a new end-to-end solution method that combines optimization and machine learning to improve solution quality while solving the tour optimization problem in 1-3 minutes. Our results from a series of computational experiments and practical case studies demonstrate consistent and substantial improvements over state-of-the-art approaches that neglect drone-flight physics. This underscores the importance of integrating accurate drone trajectory modeling into delivery planning systems to realize the full potential of truck-and-drone collaborative logistics, resulting in annual savings worth millions of dollars, while also benefiting the environment. 

Future research could extend our approach in many promising directions. One can incorporate additional practical considerations for both truck and drone—such as, variations in wind conditions, which requires adding an extra constraint in the drone model. Extending to other system configurations—for example, adding multiple drones or increasing the number of drone flights from and to the same nodes—could enhance operational flexibility. Lastly, while we focus on drone energy savings only as a secondary objective, a more holistic approach to studying the energy costs of electric trucks and advanced battery-powered drones may yield valuable insights into the overall environmental benefits and carbon reduction strategies.

\bibliographystyle{informs2014} 

\bibliography{References}

\newpage
\begin{APPENDICES}
\section{Supplemental Figures for the Computational Results}\label{section_appendix_comp_plot}

\begin{figure}[hbt]
    \centering
    \includegraphics[scale=0.4]{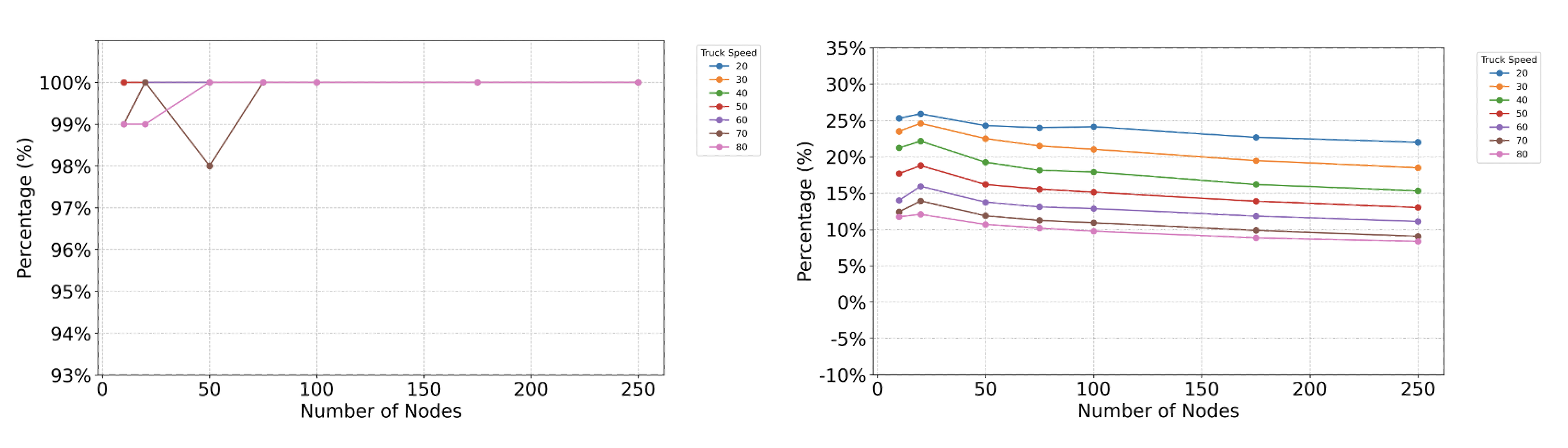}
    \caption{Duration reduction: truck-and-drone tour (P method) vs. Truck-only TSP for Scenario I (without RAS). The left sub-figure shows the percentage of instances for which the P method strictly improves upon the Truck-only TSP duration. The right sub-figure shows the average percentage reduction achieved by the P method compared with the Truck-only TSP.}
    \label{fig:s1_pvt2}
\end{figure}

\begin{figure}[hbt]
    \centering
    \includegraphics[scale=0.4]{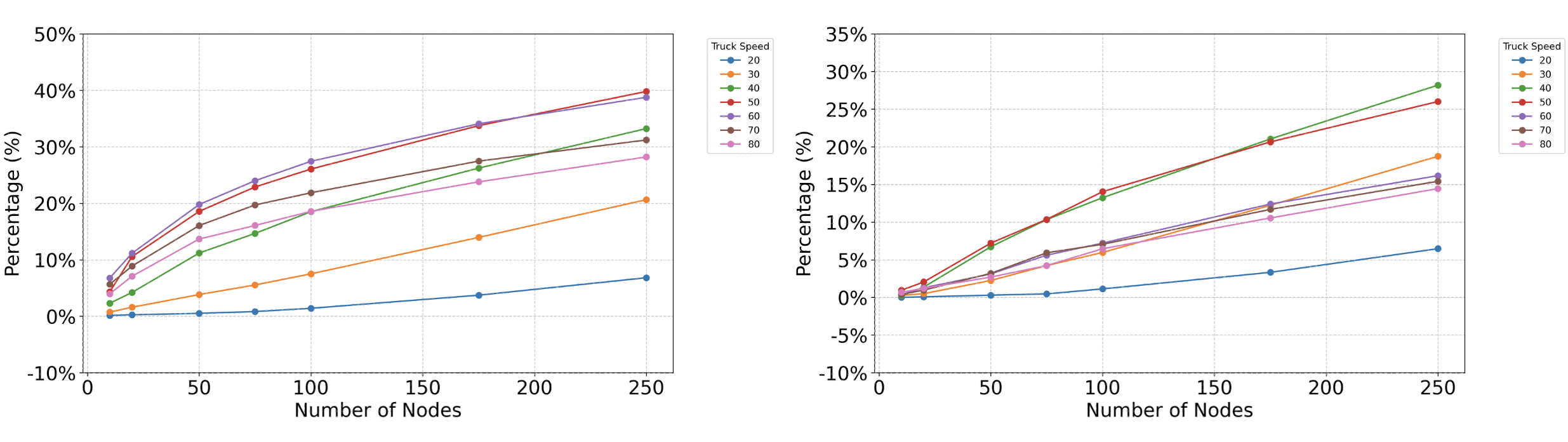}
    \caption{Average percentage reduction in the total tour duration achieved by the P method compared to the K method (left sub-figure) and MK method (right sub-figure) for Scenario I (without RAS).}
    \label{fig:s1_pvkvmk_time}
\end{figure}

\begin{figure}[hbt]
    \centering
    \includegraphics[scale=0.4]{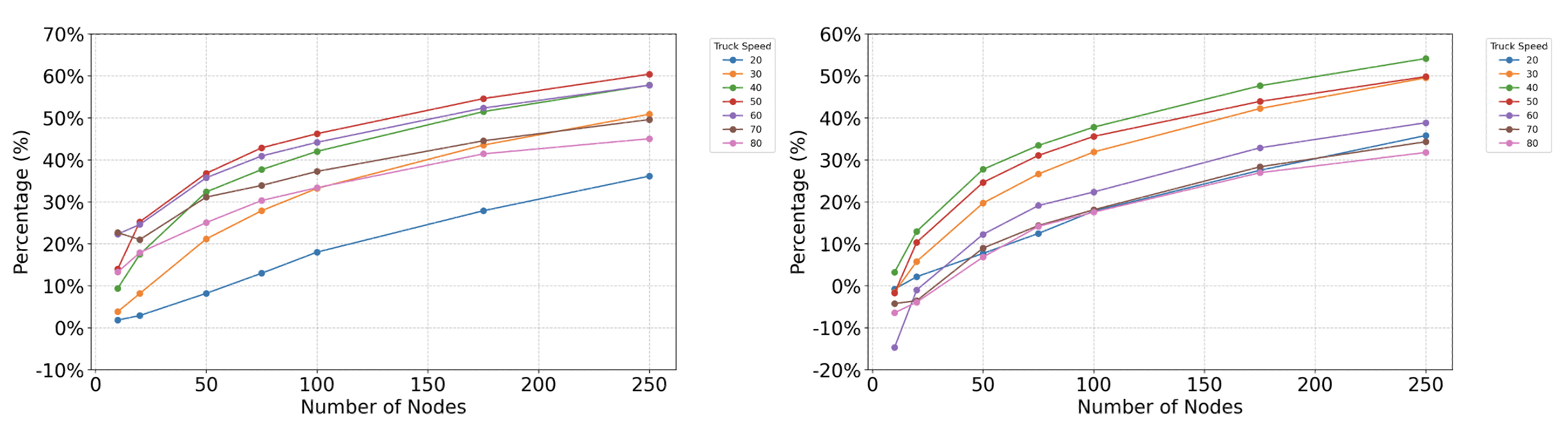}
    \caption{Average percentage reduction in drone energy consumption achieved by the P method compared to the K method (left sub-figure) and MK method (right sub-figure) for Scenario I (without RAS).}
    \label{fig:s1_pvkvmk_energy}
\end{figure}

\begin{figure}[hbt]
    \centering
    \includegraphics[scale=0.4]{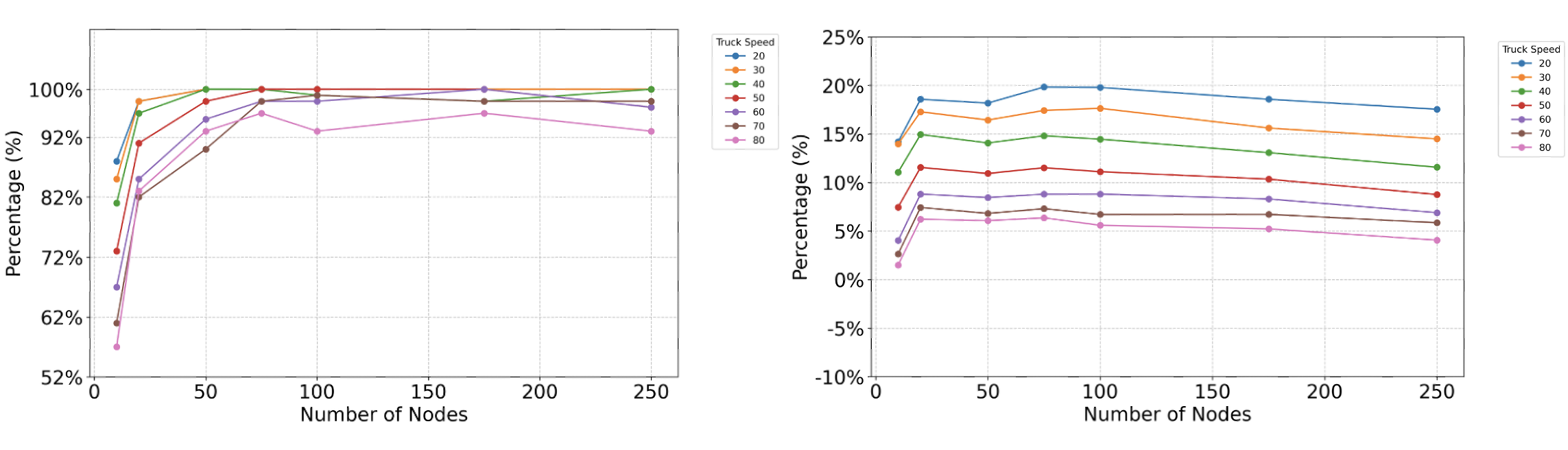}
    \caption{Duration reduction: truck-and-drone tour (P method) vs. Truck-only TSP for Scenario II (with RAS). The left sub-figure shows the percentage of instances for which the P method strictly improves upon the Truck-only TSP duration. The right sub-figure shows the average percentage reduction achieved by the P method compared with the Truck-only TSP.}
    \label{fig:s2_pvt2}
\end{figure}

\begin{figure}[hbt]
    \centering
    \includegraphics[scale=0.4]{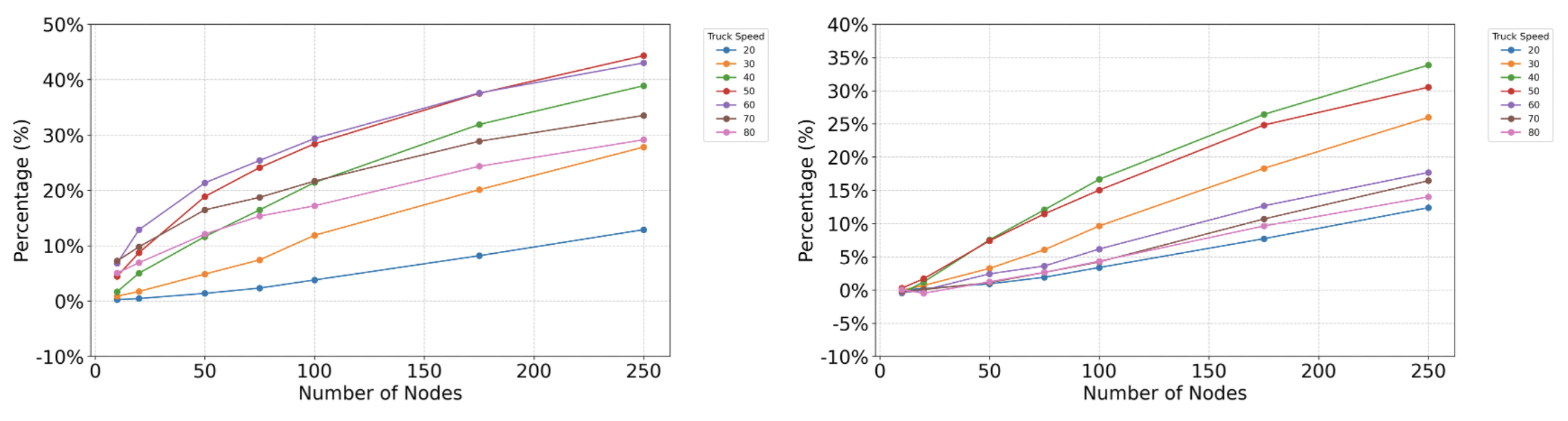}
    \caption{Average percentage reduction in the total tour duration achieved by the P method compared to the K method (left sub-figure) and MK method (right sub-figure) for Scenario II (with RAS).}
    \label{fig:s2_pvkvmk_time}
\end{figure}

\begin{figure}[hbt]
    \centering
    \includegraphics[scale=0.4]{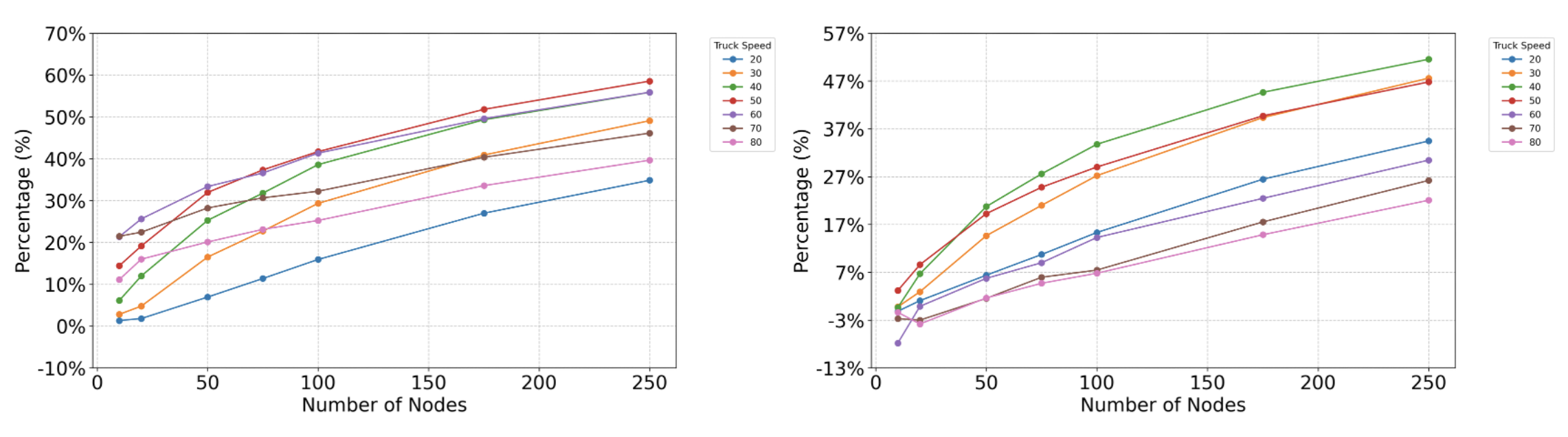}
    \caption{Average percentage reduction in drone energy consumption achieved by the P method compared to the K method (left sub-figure) and MK method (right sub-figure) for Scenario II (with RAS).}
    \label{fig:s2_pvkvmk_energy}
\end{figure}

\end{APPENDICES}


%
%
%




\end{document}